\documentclass[a4paper,12pt]{article}
\usepackage{amsmath, amssymb, amsthm, xcolor, bbm, graphicx, subcaption}
\usepackage{enumitem}
\usepackage{hyperref}
\usepackage{setspace}
\usepackage{geometry}
\usepackage{multicol}
 \usepackage[width=0.9\textwidth, margin=11mm]{caption}
\DeclareMathOperator*{\argmax}{arg\,max}

\usepackage{caption}

\usepackage[longnamesfirst]{natbib}
\bibpunct[ ]{(}{)}{,}{a}{}{,}

\newtheorem{innercustomgeneric}{\customgenericname}
\providecommand{\customgenericname}{}
\newcommand{\newcustomtheorem}[2]{%
  \newenvironment{#1}[1]
  {%
   \renewcommand\customgenericname{#2}%
   \renewcommand\theinnercustomgeneric{##1}%
   \innercustomgeneric
  }
  {\endinnercustomgeneric}
}

\newcustomtheorem{customthm}{Theorem}
\newcustomtheorem{customlemma}{Lemma}
\newcustomtheorem{customcorr}{Corollary}
\hypersetup{
    colorlinks=true,
    linkcolor=blue,
    filecolor=magenta,      
    urlcolor=cyan,
    citecolor=teal,
    pdftitle={Overleaf Example},
    pdfpagemode=FullScreen,
    }

\title{Convergence of Discrete Percolation Models to the Brownian Web Distance}
\author{Craig Belair}
\date{September 2025}
\usepackage{fancyhdr}
\begin{document}

\maketitle

\abstract{ \noindent The Brownian web is a collection of coalescing Brownian motions started from every space-time point in $\mathbb{R}^{2}$. The Brownian web can be constructed as a scaling limit of coalescing one-dimensional simple random walks started at every point in a two-dimensional space-time lattice. \cite{VV} introduced a family of discrete random distance functions defined on these sequences of rescaled lattices. It was shown that, given the appropriate notion of convergence, these discrete distance functions converge to a function known as the Brownian web distance. We introduce a new method of argument that allows us to show that a broad class of discrete first passage percolation models also converge to the Brownian web distance. Unlike the arguments used in \cite{VV}, our methods do not depend on the use of planar dual graphs. This allows our methods to be applied to models that allow random walks to cross one another before coalescing.}
\section{Introduction}
\subsection{Background}
The Brownian web is a random network consisting of coalescing one-dimensional Brownian motion paths starting from every space-time point in $\mathbb{R}^{2}$. Based on the work of \cite{AR}, the Brownian web was first rigorously constructed by \cite{TW}. \cite{FINR} later developed convergence criteria for the Brownian web which is analogous to Donsker's invariance principle (\cite{DON}.) In particular, they showed that the Brownian web can be constructed as the scaling limit of a web of coalescing simple random walks started at every point in a two-dimensional space-time lattice. \cite{NRS} extended these results by showing that, under certain conditions, random walks that allow crossing before coalescence exhibit the same convergence.
\newline
\newline
A directed first passage percolation model, $D_{BR}$, can be naturally defined for the Brownian web. In this model, a particle can freely follow the forward time trajectory of a path in the Brownian web, or can jump to a different path started at the particle's current space-time point. For $u,v \in \mathbb{R}^{2}$, the distance $D_{BR}(u, v)$ is defined as the minimum number of jumps required to travel from $u$ to $v$. Importantly, there exist special random space-time points at which multiple (at most three) paths in the Brownian web originate (shown in \cite{FINR} among other places.) This ensures that the model is non-trivial, i.e. that jumps can occur. In this paper, we examine analogous models defined on discrete webs. We show that under appropriate rescaling, our discrete directed first passage percolation models weakly converge to $D_{BR}$.
\subsection{Our contributions}
We begin our paper by defining a directed first passage percolation model using the webs of non-simple coalescing random walks introduced in \cite{NRS}. We denote this percolation model $D_{RW}$. The random walk webs used to construct our model differ from those used in \cite{VV}, as we allow paths to cross each other before coalescing. In our model, we use a deterministic set $\mathcal{K} \subset \mathbb{Z}^{2}$ to capture the `jumps' that a particle is permitted to make. A particle can freely follow the forward-time trajectory of any walk in the web. The particle may also jump from $x \in \mathbb{Z}^{2}$ to $y \in \mathbb{Z}^{2}$ if and only if $y-x \in \mathcal{K}$. For $u,v \in \mathbb{Z}^{2}$, $D_{RW}(u, v)$ is defined as the minimal number of jumps required for a particle to travel from $u$ to $v$. We define $D_{RW}^{n}$ as the rescaled distance functions, with space rescaled by a factor of $\sqrt{cn}$ and time rescaled by a factor of $n$, where $c$ is the variance of the random walk increment. In this setting, we show the following:

\begin{customthm}{1.1}\label{T1.1}
Let $\mathcal{V}^{n}$ be a system of coalescing random walks started at every point in $\mathbb{Z}^{2}$, having increment $\zeta$, and such that the walks are independent until meeting, at which point they coalesce. Suppose further that $\zeta$ satisfies:
\begin{itemize}
\item $\mathbb{E}(|\zeta|^{5}) < \infty$
\item $\mathbb{E}(\zeta) = 0$
\item A random walk with increment $\zeta$ is aperiodic
\end{itemize}
Then there exists a coupling of $\mathcal{V}^{n}$ and the Brownian web, $\mathcal{W}$, such that $D_{RW}^{n}$ converges (in the epigraph sense) to $D_{BR}$ as $n \to \infty$.
\end{customthm}
\noindent Epigraph convergence is a standard means of characterizing the convergence of functions which are lower semicontinuous but not continuous. A detailed definition of epigraph convergence is given in \hyperref[S2.1]{Section 2.1}. We deduce Theorem \ref{T1.1} by showing that the hypotheses of Lemma \ref{4.1VV} hold. The reader may be interested in referring to these hypotheses, as they provide an alternative way of understanding the results of Theorem \ref{T1.1}.
\newline
\newline
\noindent A priori this epigraph convergence is not obvious. The issue is that microscopic information, such as two curves being only a few vertices apart from one another, does not behave well in the scaling limit (which only sees two curves being less-than-scaling away from each other.) For an example, see \hyperref[Fig1]{Figure 1}.
\begin{figure}[htp]    
    \centering
    \includegraphics[width=8cm]{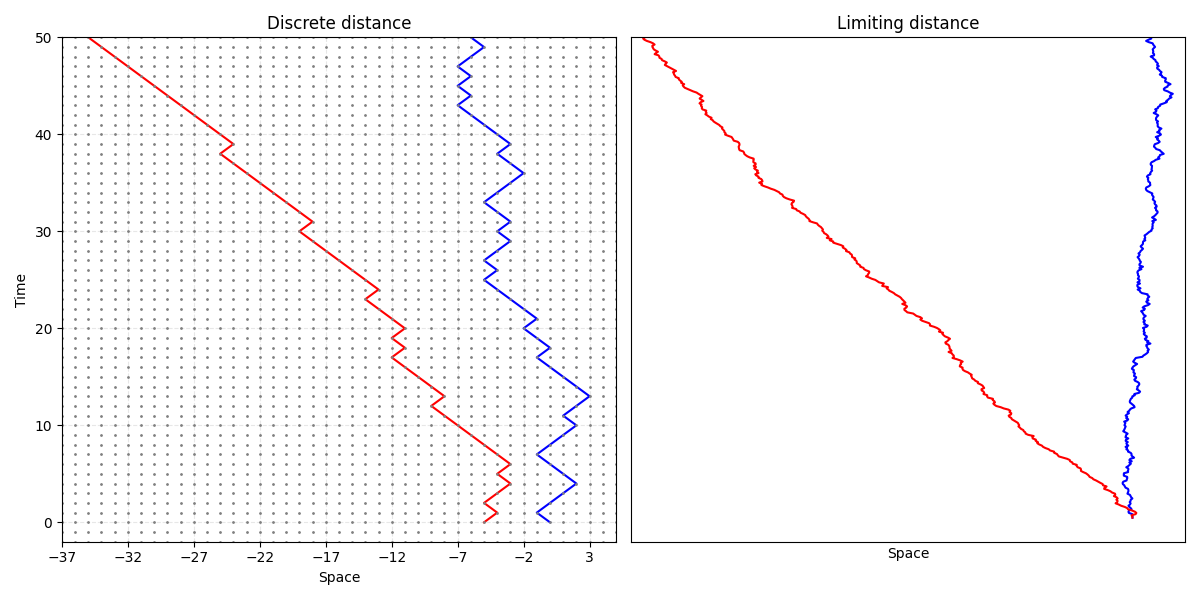}
    \caption{The discrete paths on the left start five vertices apart from one another. In a discrete model which allows a particle to jump only to directly adjacent vertices, a particle cannot travel from the red curve to the blue curve using a single jump. However, since $5 \ll \sqrt{n}$, the distance between the start points of red and blue paths can disappear in the scaling limit. As such, in the scaling limit, the red and blue paths may start from the same space-time point, hence a particle can jump between the paths in the model $D_{BR}$. }
\end{figure}\label{Fig1}
\text{ }
\newline
\newline
Our approach is to introduce the concepts of journeys and itineraries. Journeys, $J(t)$, are constructed by first selecting a point $(x, s) \in \mathbb{Q}^{2}$ and following along the path in a web (either the discrete web or the Brownian web) started at space-time point $(x,s)$ (or at a nearby lattice point, in the discrete case). For rational times $\sigma_{i}$ satisfying $s < \sigma_{1} < \sigma_{2} < ... < \sigma_{k}$ we choose to follow either the right or left forward-time boundary curve of the point $(J(\sigma_{i}), \sigma_{i})$. We follow that boundary curve up until time $\sigma_{i+1}$ and then once again choose to follow either the left or right boundary curve from that point. An example of a journey on the Brownian web is depicted in \hyperref[Fig6]{Figure 6} in the definitions section. Journeys themselves are used to approximate geodesics, which allows us to focus on points lying on these boundary curves, rather than points lying between the left and right boundary curves of a point or spatial interval. This avoids the issue relating to forward-time boundary curves which is described in the caption of \hyperref[Fig3]{Figure 3}. The use of rational start points and times ensures that the set of journeys is countable. An itinerary is the list of rational numbers which parameterize a journey. Itineraries provide convenient language for pairing discrete journeys with their counterparts in the Brownian web.
\newline
\newline
We first prove a type of modified pointwise convergence of discrete journeys to the corresponding journeys on the Brownian web. We do this by introducing approximation curves, which behave similarly to the discrete journey, but are easier to analyze. We then show that journeys defined on the Brownian web can provide arbitrarily good approximations for $D_{BR}$ geodesics. Then, using Lemma \ref{4.1VV}, we combine these two results to show that the epigraphs of $D_{RW}^{n}$ converge to those of $D_{BR}$. Lastly, we extend the pointwise convergence of discrete journeys to a stronger type of convergence -- uniform convergence on compact sets. We show this stronger convergence by analyzing a modulus of continuity for our approximation curves, and then using it to recover a modulus of continuity for our discrete journeys.
\newline
\newline
The main results of our paper are specific to models constructed on sequences of rescaled lattices; however the method of using journeys and itineraries is quite general and may be applied to many other percolation models. Journeys provide a way of using properties of the forward-time boundary curves to recover important information about a first passage percolation model. In Theorem \ref{T6.1.1} we show that any $D_{BR}$ geodesic can be well approximated by a journey, thereby showing that the behavior of $D_{BR}$ can be effectively understood through the behavior of Brownian web journeys.
\subsection{Additional remarks}
We do not believe the conditions that we impose on the set $\mathcal{K}$ are optimal. In our proof, the require $\mathcal{K}$ be a bounded set satisfying:
\begin{enumerate}
\item $(x,0) \in \mathcal{K}$ for some $x<0$
\item $(x,0) \in \mathcal{K}$ for some $x>0$
\item $t \geq 0$ for all $(x,t) \in \mathcal{K}$
\end{enumerate}
However, we suspect that our results can be achieved by requiring only that $\mathcal{K}$ is bounded and contains at least one point $(x,t)$ such that $(x,t) \neq (0,0)$. 
\newline
\newline
If the random walk increment $\zeta$ described in Theorem \ref{T1.1} is assumed to be bounded, then our results can be used to show that models allowing for $(x,t) \in \mathcal{K}$ with $t<0$ exhibit the convergence given in Theorem \ref{T1.1}. This extension can be achieved by considering a jump set $\mathcal{K}$ containing negative time jumps, and comparing it to a strategically selected set $\mathcal{K}^{\star}$, containing only non-negative time jumps. Due to the assumed boundedness of $\zeta$, one can choose $\mathcal{K}^{\star}$ as so that the any discrete journey determined using the jump set $\mathcal{K}$ is dominated by the corresponding journey determined using the jump set $\mathcal{K}^{\star}$. This domination can then be used to show that the convergence of the first passage percolation model corresponding to $\mathcal{K}^{\star}$ implies convergence of the model corresponding to $\mathcal{K}$.
\subsection{Related work}
In \cite{FINR}, the Brownian web is constructed as the scaling limit of discrete webs of coalescing simple random walks. When simple walks are used, the graphs of the discrete webs are planar. This planarity allows for a naturally defined dual web, consisting of backwards walks (see \hyperref[Fig2]{Figure 2}). Up to a spatial shift and time reversal, the discrete dual web has the same distribution as the original discrete web. Further, any realization of the discrete web induces a unique realization of the dual web. It turns out this duality is preserved in the limit: the Brownian web also has a dual web. Both the discrete dual web and the Brownian web dual have the special property that no path in the primal can cross a path in the dual. These dual webs have proven instrumental in understanding properties of the Brownian web and developing convergence criteria, for example in \cite{FINR}, \cite{NRS}, and \cite{ESS}.
\begin{figure}[htp]    
    \centering
    \includegraphics[width=8cm]{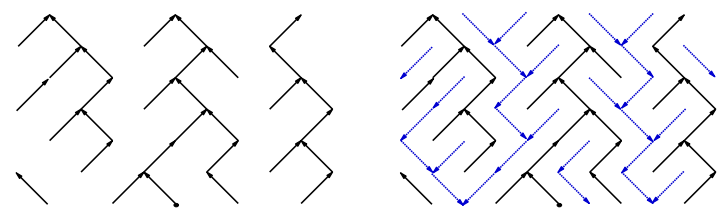}
    \caption{A discrete web defined on the even integer lattice shown in black, and the associated dual web, defined on the odd integer lattice, shown in blue. Figure taken from \cite{ESS}}
\end{figure}\label{Fig2}
\text{ }
\newline
\newline
\noindent \cite{VV} defined a first passage percolation model on the discrete webs arising from simple random walks. They showed epigraph convergence of their models to the Brownian web distance. The method used in \cite{VV} was to analyze the backward-time boundary curves from a spatial interval $I$, and to compare those curves to the discrete analogs. Through arguments using the discrete dual web and the Brownian web dual, it was shown that the discrete boundary curves converged to those in the Brownian web. The use of backward-time boundary curves was not a stylistic choice -- the information given by the backwards boundaries differs from the forward boundaries (see \hyperref[Fig3]{Figure 3}.)
\newline
\newline
The bulk of our arguments necessarily differ from the arguments used in \cite{VV}. The main issue is that the arguments used in \cite{VV} rely heavily on the analysis of the discrete dual web. When non-simple random walks are considered, the random walk paths may cross many times before coalescing. As a result, the graph of these discrete webs is no longer planar, and there is no natural definition of a discrete dual web (see \hyperref[Fig4]{Figure 4}). This makes the methods used in \cite{VV}, which compared the discrete dual to the Brownian dual, inapplicable.
\newline
\newline
In addition to the work done in \cite{VV}, our research is motivated in part by the results of \cite{NRS}. The work of \cite{NRS} showed that, provided certain conditions were met, non-simple random walk webs also converge to the Brownian web. In particular, \cite{NRS}, provided a convergence criteria for discrete random walk webs which allow for paths to cross one another. \cite{NRS} showed the following (Theorem 1.3.3 in \cite{NRS}):
\newline
\newline
\noindent Let $\mathcal{V} := \{ Y_{x,t} \}$ be a collection of coalescing random walks started from each point $(x,t) \in \mathbb{Z}^{2}$, such that the walks are independent before they meet, coalesce upon meeting, and have increments distributed as $\zeta$. Suppose $\zeta$ satisfies:
\begin{itemize}
\item $\mathbb{E}(|\zeta|^{5}) < \infty$
\item $\mathbb{E}(\zeta^{2})=c$
\item $\mathbb{E}(\zeta) = 0$
\item A random walk with increment $\zeta$ is aperiodic
\end{itemize} 
Then, under $\sqrt{cn}, n$ space-time rescaling, the rescaled webs $\mathcal{V}^{n}$ will converge weakly to the Brownian web. Here convergence takes place in the space $(\mathcal{H}, d_{\mathcal{H}})$, described in \hyperref[S2.1]{Section 2.1}.
\newline
\newline
The results of \cite{NRS} and those of \cite{VV} may lead one to believe that first passage percolation models defined on discrete webs that allow for crossing paths may also exhibit the convergence shown in \cite{VV}. In this paper we prove this intuition to be correct.

\begin{figure}[htp]    
    \centering
    \includegraphics[width=8cm]{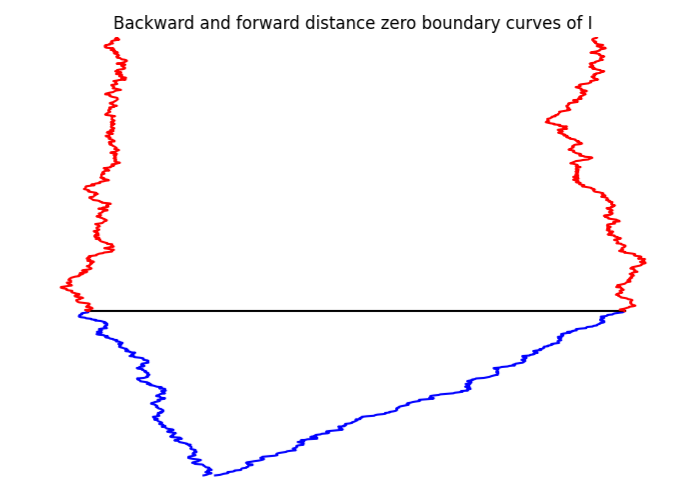}
    \caption{A spatial interval $I$ is depicted in black. The backward-time distance zero boundary curves from I are depicted in blue, while the forward-time distance zero boundary curves are depicted in red. Any point $u$ between the blue curves always satisfies $D_{BR}(u, v) = 0$ for some $v \in I$. This is because at any such point $u$, there originates a Brownian motion path and that path cannot cross the blue boundary curves, lest it coalesce with them. However, for many (in fact, almost all) points $v$ between the red curves, there are no points $u \in I$ such that $D_{BR}(u, v) = 0$. This is because the paths originating along $I$ rapidly coalesce with one another, so they cover very little space. The same idea holds for higher distance boundaries.}
\end{figure}\label{Fig3}
\begin{figure}[h! tbp]
  \centering
  \begin{subfigure}[b]{0.4\linewidth}
    \includegraphics[width=\linewidth]{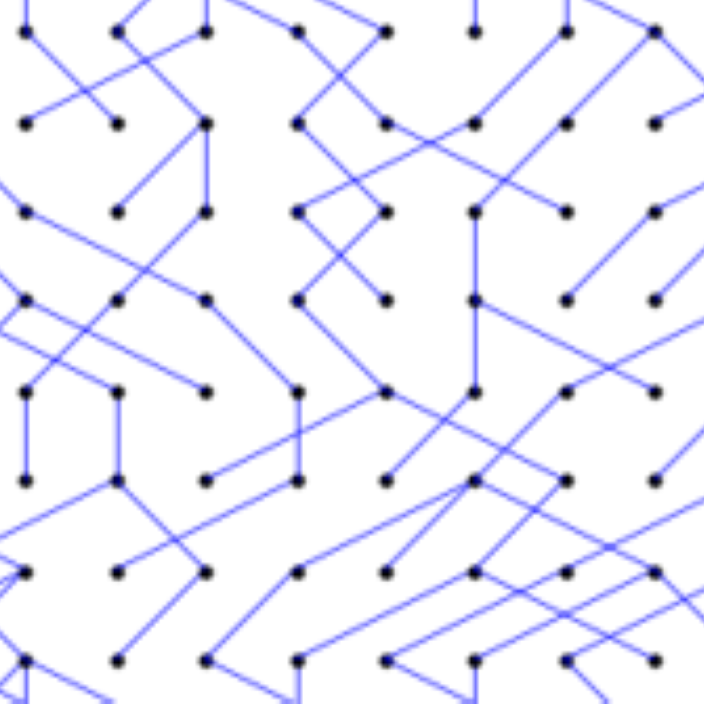}
    \caption{A discrete random walk web consisting of non-simple coalescing random walks. Unlike in the case of \hyperref[Fig2]{simple random walks}, there is no natural way to define a dual web. This is because walkers may cross one another.}
  \end{subfigure}
  \begin{subfigure}[b]{0.4\linewidth}
    \includegraphics[width=\linewidth]{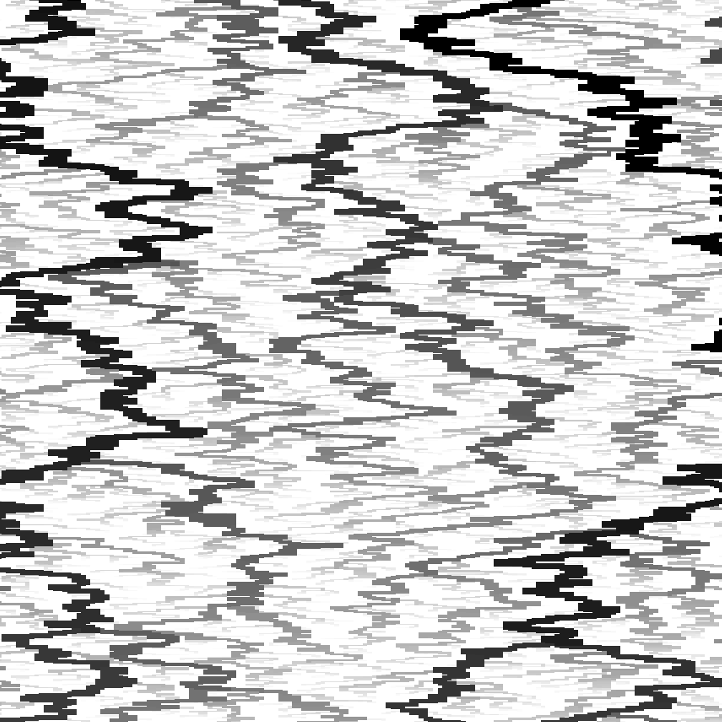}
    \caption{An approximation of a realization of the Brownian Web. As shown by \cite{NRS}, discrete webs consisting of non-simple walks may still converge to the Brownian web. In this case, the crossing of paths disappears in the limit.}
  \end{subfigure}
  \caption{}
\end{figure} \label{Fig4}
\newpage
\section{Definitions and preliminaries}
\subsection{Basic definitions} \label{S2.1}
\label{conditions}
Let $\zeta$ be a centered, integer-valued random variable with finite positive variance. $\zeta$ will serve as the increment variable for a family of random walks started at every point on $\mathbb{Z}^{2}$. To avoid complicated notation, throughout the paper we assume $\zeta$ has unit variance. By replacing the spatial rescaling factor of $\sqrt{n}$ by $\sqrt{cn}$ where $c$ is the variance of $\zeta$, one can obtain the same results for $\zeta$ with arbitrary finite positive variance.
\newline
\newline
Let $\{ \zeta_{i, j} \}_{(i, j) \in \mathbb{Z}^{2}}$ be a family of i.i.d. random variables distributed as $\zeta$. We use this family to define a random walk started at every point in $\mathbb{Z}^{2}$. The first coordinate represents space, while the second represents time. For $(i, j) \in \mathbb{Z}^{2}$ we let $Y_{i, j} : [j, \infty) \cap \mathbb{Z} \to \mathbb{Z}$ be the function defined as follows:
$$Y_{i,j}(t) := \begin{cases} i & \text{if } t=j \\ Y_{i+\zeta_{i,j}, j+1}(t) & \text{else} \end{cases}$$
\newline
We extend the function by performing a linear interpolation between integer time points and setting the function equal to zero at infinity. That is:
$$Y_{i, j}(a) := \begin{cases} (t+1-a)Y_{i, j}(t) + (a-t)Y_{i, j}(t+1) & \text{for } t \in \mathbb{Z} \cap [j, \infty), a \in [t, t+1) \\ 0 & \text{for } a = \infty \end{cases}$$
\newline
Next, to examine the scaling limits, we define the rescaled random walk functions. For $(\sqrt{n}i, nj) \in \mathbb{Z}^{2}$, we define functions $Y^{n}_{i, j} : [j, \infty] \to \mathbb{R}$ as:
$$Y^{n}_{i, j}(t) := \frac{1}{\sqrt{n}}Y_{\sqrt{n}i, nj}(nt)$$
For technical reasons, it is useful to introduce functions $\gamma_{\pm \infty, r} : [r, \infty] \to \overline{\mathbb{R}}$ as $$\gamma_{\pm \infty, r}(t) := \pm \infty$$ and define a set $\mathcal{P}^{n} := \{ s : ns \in \mathbb{Z} \text{ or } s = \pm \infty \}$. With these definitions in hand, we give the law of random variables $\mathcal{V}^{n}$, whose distribution is equal to that of the collection of random walk paths on a rescaled lattice:
$$\mathcal{V}^{n} \overset{d}{:=} \{ Y^{n}_{i, j} : (\sqrt{n}i, nj) \in \mathbb{Z}^{2} \} \cup \{ \gamma_{\infty,r}, \gamma_{- \infty,r} :r \in \mathcal{P}^{n} \}$$
It can be seen that $\{ Y^{n}_{i, j} \}$ form a collection of coalescing one-dimensional random walks. The inclusion of the paths $\gamma_{\pm \infty, r}$ will allow us to define $\mathcal{V}^{n}$ as a random variable taking values in a particular metric space, which will later be defined.
\newline
\newline
The Brownian web, which we denote $\mathcal{W}$, describes a collection of one-dimensional coalescing Brownian motions started from every space-time point in $\mathbb{R}^{2}$. The paths in $\mathcal{W}$ evolve independently of one another until the point at which they first meet. Upon meeting, paths will coalesce with one another. A detailed characterization of the Brownian web is given in \cite{FINR}. $\mathcal{V}^{n}$ and $\mathcal{W}$ can both be viewed as random variables taking values in a metric space $(\mathcal{H}, d_{\mathcal{H}})$. In an effort to describe this metric space, we first give some intermediary definitions. These definitions were previously given in \cite{FINR}.
\newline
\newline
Let $\Phi : \overline{\mathbb{R}}^{2} \to [-1,1]$ be defined as:
$$\Phi( x,t) := \begin{cases} \frac{\tanh(x)}{1+|t|} & \text{for } t \in \mathbb{R} \\ 0 & \text{for } t \in \{ \pm \infty \} \end{cases}$$
\newline
With the convention that $\tanh(\pm \infty) = \pm 1$.
\newline
\newline
For $t_{0} \in [-\infty, \infty]$ let $C[t_{0}]$ denote the set of functions $f : [t_{0}, \infty] \to \overline{\mathbb{R}}$ such that $\Phi(f(t),t)$ is continuous. Define the space $\Pi$ as:
$$\Pi := \bigcup\limits_{t_{0} \in [-\infty, \infty]} C[t_{0}] \times \{ t_{0} \}$$
For $(f, t_{0}) \in \Pi$ we let $\hat{f}$ represent the function obtained by setting $\hat{f}(t) = f(t_{0})$ for $t \leq t_{0}$. We define a distance function on $\Pi$ as follows:
$$d((f_{1}, t_{1}), (f_{2}, t_{2})) := \sup\limits_{t \in \overline{\mathbb{R}}} \left| \Phi(\hat{f_{1}}(t),t) - \Phi(\hat{f_{2}}(t),t) \right| \vee |\tanh(t_{1}) - \tanh(t_{2})|$$
We afford ourselves common abuses of notation relating to the space $(\Pi, d)$. In the case where $t_{0}$ is well understood or not pertinent, we may use the shorthand $f \in \Pi$ to denote the path $(f, t_{0}) \in \Pi$. Further, for functions $f :[t_{0}, \infty) \to \overline{\mathbb{R}}$ we may consider these functions as elements of $\Pi$. In such cases, we are implicitly extending $f$ by defining $f(\infty) = 0$.
\newline
\newline
Note that by definition of the space $(\Pi, d)$, two paths $f_{1}, f_{2} \in \Pi$ may be identified even if it is the case that $f_{1}(\pm \infty) \neq f_{2}( \pm \infty)$.
\newline
\newline
We use $\Pi^{\star} \subset \Pi$ to denote the set of all paths with finite starting times, taking only finite spatial values at all finite times. That is, paths $(f, t_{0}) \in \Pi$ such that $t_{0} \in \mathbb{R}$ and $\hat{f}(t) \neq \pm \infty$ for all $t \in \mathbb{R}$. Due to the assumed continuity of $\Phi(f(t),t)$, if $f \in \Pi^{\star}$ then the restriction of $\hat{f}$ to $\mathbb{R}$ is continuous with respect to the standard topology on $\mathbb{R}$.
\newline
\newline
$(\Pi, d)$ is a complete separable metric space, as shown in \cite{FINR}. We let $\mathcal{H}$ denote the collection of compact subsets of $(\Pi, d)$, and equip $\mathcal{H}$ with the Hausdorff metric, $d_{\mathcal{H}}$:
$$d_{\mathcal{H}} (K_{1}, K_{2}) := \sup\limits_{g_{1} \in K_{1}} \inf\limits_{g_{2} \in K_{2}} d(g_{1}, g_{2}) \vee \sup\limits_{g_{2} \in K_{2}} \inf\limits_{g_{1} \in K_{1}} d(g_{1}, g_{2})$$
$(\mathcal{H}, d_{\mathcal{H}})$ is also a complete separable metric space (also shown in \cite{FINR}). 
\newline
\newline
For a function $h : \mathbb{R}^{4} \to \overline{\mathbb{R}}$ we denote the epigraph of $h$ as $\mathfrak{e}(h)$:
$$\mathfrak{e}(h) := \{ (x,y) \in \mathbb{R}^{4} \times \overline{\mathbb{R}} : y \geq h(x) \}$$
We define a map $E: \mathbb{R}^{4} \times \overline{\mathbb{R}} \to \mathbb{R}^{4} \times [-1, 1]$ as:
$$E(r,s) := \begin{cases} \left( r, \frac{|u-v|e^{-|r|}s}{1+|s|} \right) & s \neq \pm \infty \\ (r, \pm|u-v|e^{-|r|}) & s = \pm \infty \end{cases}$$
Where $r=(u,v)$ with $u, v \in \mathbb{R}^{2}$.
\newline
\newline
We equip the space $\mathbb{R}^{4} \times \overline{\mathbb{R}}$ with the topology induced by the map $E$. Under this topology, we consider a collection of subsets, $\mathcal{E}_{\star}$:
$$\mathcal{E}_{\star} := \{ \Gamma \subseteq \mathbb{R}^{4} \times \overline{\mathbb{R}} : \Gamma \text{ is closed}, \Gamma \cap (\{x \} \times \overline{\mathbb{R}}) \neq \emptyset \text{ }\forall x \in \mathbb{R}^{4} \}$$
For $\Gamma_{1}, \Gamma_{2} \in \mathcal{E}_{\star}$ we define $d_{\star}(\Gamma_{1}, \Gamma_{2})$ to be the Hausdorff distance of $E(\Gamma_{1}), E(\Gamma_{2})$. The space $(\mathcal{E}_{\star}, d_{\star})$ is compact by Lemma 7.1 in \cite{DV}.
\newline
\newline
$(\mathcal{E}_{\star}, d_{\star})$ serves as a suitable space to examine the convergence of the epigraphs of our percolation models. We now define these percolation models.
\subsection{Percolation models}
Let $\mathcal{K}$ be any deterministic bounded subset of $\mathbb{Z}^{2}$ satisfying the following:
\begin{itemize}
\item $(x, 0) \in \mathcal{K}$ for some $x > 0$
\item $(x, 0) \in \mathcal{K}$ for some $x<0$
\item $t \geq 0$ for all $(x, t) \in \mathcal{K}$
\end{itemize}
 $\mathcal{K}$ will be used to determine the `open edges' in our discrete first passage percolation models. That is, the set of vertices that the model allows for jumps between; if $v \in \mathcal{K}$, then a particle may travel from $(0, 0)$ to $v$ at a cost of at most one. For $u = (u_{1}, u_{2}) \in \mathbb{Z}^{2}$, define the translated set:
$$\mathcal{K}_{u} := \{ (x,t) : (x-u_{1},t-u_{2}) \in \mathcal{K} \}$$
In order to explicitly define our weight function, we introduce two auxiliary variables. For $(i, j), (\ell, k) \in \mathbb{Z}^{2}$, define:
\newline
$p_{(i, j), (\ell, k)} := \begin{cases} 0 & \text{if } (\ell, k) = (i, j) \\ 0 & \text{if } (\ell, k) = (i + \zeta_{i, j}, j+1) \\ \infty & \text{else} \end{cases}$
\newline
$s_{(i, j), (\ell, k)} := \begin{cases} 1 & \text{if } (\ell, k) \in \mathcal{K}_{(i, j)} \\ \infty & \text{else} \end{cases}$
\newline
These definitions allow us to succinctly describe the weight function, $w_{u, v}$. For $u, v \in \mathbb{Z}^{2}$, the weight function $w_{u, v}$ denotes the weight of the directed edge connecting $u \to v$:
$$w_{u,v} := p_{u, v} \wedge s_{u, v}$$
\begin{figure}[htp]
    \centering
    \includegraphics[width=8cm]{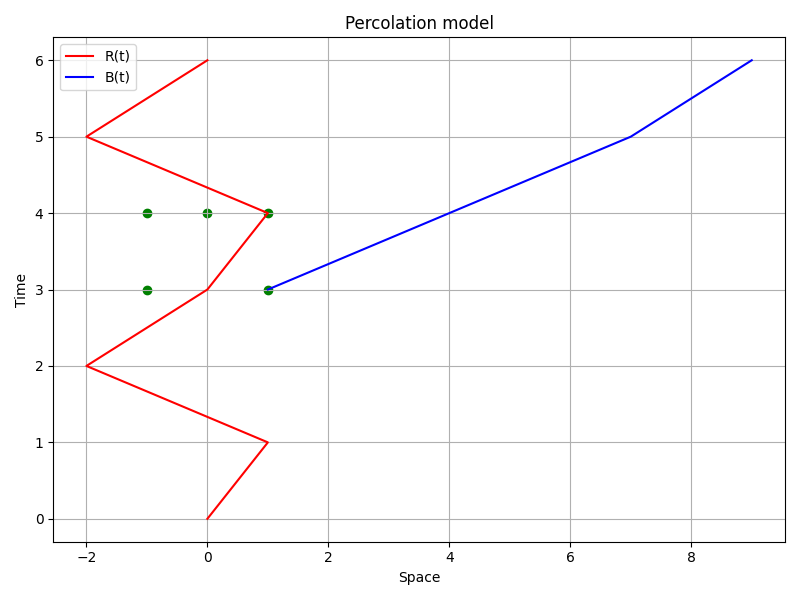}
    \caption{A depiction of the percolation model on two random walks. $\mathcal{K}$ is chosen to be the five vertices adjacent to and `above' $(0,0)$:
    $\mathcal{K} = \{ (-1, 0), (1,0), (-1,1), (0,1), (1,1) \}$.
      Since $(1, 3) \in \mathcal{K}_{(0,3)}$, a jump can be made from $R(t)$ to $B(t)$ at the time $t=3$ for a cost of one. As a result, the first passage percolation distance from $(0,0) \to (9,6)$ is one. The first passage percolation distance from $(0,0) \to (0,6)$ is zero, as a walker can freely follow the upward path of $R(t)$. Note that we can only follow the upward path freely. As a result, the model is asymmetric in time; the first passage percolation distance from $(0,6) \to (0,0)$ would instead be infinity.}
    \label{fig:pmodel}
\end{figure} 
We then define the first passage percolation distance. For $u,v \in \mathbb{R}^{2}$, and random variable $\mathcal{V}^{1}$ with the associated increment variables $\zeta_{i, j}$, we define:
$$D_{RW}^{1}(u, v) := \begin{cases} \inf \left\{ \sum_{i=0}^{i=m}w_{x_{i}, x_{i+1}} : m \in \mathbb{Z}_{\geq 0}, x_{0}=u, x_{m+1} = v \right\} & u, v \in \mathbb{Z}^{2} \\ \infty & \text{else} \end{cases}$$
Using this definition, $D_{RW}^{1}$ is a random function which is uniquely determined by $\mathcal{V}^{1}$. We can similarly define $D_{RW}^{n}$ as a random function determined by the rescaled lattice $\mathcal{V}^{n}$. That is, the law of $D_{RW}^{n}$ is defined as follows:
$$D_{RW}^{n}((i, j), (\ell, k)) \overset{d}{=} D_{RW}^{1}((\sqrt{n}i, nj), (\sqrt{n} \ell, nk))$$
The conjectured limit function, the \textit{Brownian web distance}, is a random function from $\mathbb{R}^{4} \to \mathbb{Z}_{\geq 0} \cup \{ \infty \}$. We call this distance function $D_{BR}$. We repeat the definition previously given in \cite{VV}:
\newline
\newline
For any $((x,s),(y,t)) \in \mathbb{R}^{4}$ we let $D_{BR}((x,s),(y,t))$ denote the infimum of non-negative integers $k$ for which there exists points $(x_{1}, t_{1}), ..., (x_{k}, t_{k}) \in \mathbb{R}^{2}$ such that $s < t_{1} < ... < t_{k} < t$ and there exists a continuous path $\pi : [s,t] \to \mathbb{R}$ so that $\pi$ restricted to the intervals $[s, t_{1}], [t_{1}, t_{2}], ..., [t_{k-1}, t_{k}], [t_{k}, t]$ coincides with trajectories in the Brownian web $\mathcal{W}$. The infimum is equal to $\infty$ if there is no such $k \in \mathbb{Z}_{\geq 0}$.
\subsection{Rational itineraries and journeys}
In this section we introduce itineraries and journeys. Journeys are functions intended to approximate geodesics in our percolation models. An itinerary is a sequence of rational numbers which describe the way a journey evolves over time. For $k \in \mathbb{Z}_{\geq 0}$, we say that $\Xi_{k}$ is a length $k$ itinerary if $\Xi_{k} = (x,s,\overline{\sigma}, \overline{\eta})$ with:
\begin{itemize}
\begin{multicols}{2}
\item $x,s \in \mathbb{Q}$
\item $\overline{\sigma} \in \mathbb{Q}^{k}$
\item $\overline{\eta} \in \{ -1, 1 \}^{k}$
\item $s < \sigma_{1} < \sigma_{2} < ... < \sigma_{k}$
\end{multicols}
\end{itemize}
Let $\hat{\Xi}_{j} = (\hat{x}, \hat{s}, \hat{\overline{\sigma}}, \hat{\overline{\eta}})$ and $\Xi_{k} = (x,s,\overline{\sigma}, \overline{\eta})$ be itineraries. We say that $\hat{\Xi}_{j}$ is a continuation of $\Xi_{k}$ (and that $\Xi_{k}$ is a restriction of $\hat{\Xi}_{j}$) if the following hold:
\begin{itemize}
\begin{multicols}{2}
\item $\hat{x}=x$
\item $\hat{s}=s$
\item $j \geq k$
\item For $1 \leq i \leq k$, $\hat{\sigma}_{i} = \sigma_{i}$ and $\hat{\eta}_{i} = \eta_{i}$
\end{multicols}
\end{itemize}
We will now use itineraries to define journeys. We first define journeys of length zero. Let $\Xi_{0} = (x,s)$ be a length zero itinerary, define $\mathcal{G}_{\Xi_{0}} : [s, \infty) \to \mathbb{R}$ as:
$$\mathcal{G}_{\Xi_{0}}(t) := \gamma_{x, s}(t)$$
Here $\gamma_{x, s} \in \mathcal{W}$ is the path in $\mathcal{W}$ starting from the space-time point $(x,s)$. As shown by Theorem 2.11 in \cite{ESS}, almost surely it is the case that for every $(x,s) \in \mathbb{Q}^{2}$ there is a \textbf{unique} path $\gamma_{x, s} \in \mathcal{W}$ started at the space-time point $(x,s)$. It follows that almost surely it is the case that $\mathcal{G}_{\Xi_{0}}$ is well defined for all rational pairs $(x,s)$. We similarly define a journey with itinerary $\Xi_{0}$ in the discrete setting. Throughout the paper we use the notation $\lceil s \rceil_{n}$ in the following way:
$$\lceil s \rceil_{n} := \frac{\lceil ns \rceil}{n}$$
$\lfloor s \rfloor_{n}$ is used similarly. We define the function $G^{n}_{\Xi_{0}} : [\lceil s \rceil_{n}, \infty) \to \mathbb{R}$ as:
$$G^{n}_{\Xi_{0}}(t) := Y^{n}_{\lceil x \rceil_{\sqrt{n}}, \lceil s \rceil_{n}}(t)$$
For $k \geq 1$, we will recursively define $\mathcal{G}_{\Xi_{k}}$ and $G^{n}_{\Xi_{k}}$. Let $\Xi_{k} = (x,s,\overline{\sigma}, \overline{\eta})$ be a length $k$ itinerary. Let $\hat{\Xi}_{k-1}$ be the $k-1$ length restriction of $\Xi_{k}$. Then define $\mathcal{G}_{\Xi_{k}} : [s, \infty) \to \mathbb{R}$ as follows:
$$\mathcal{G}_{\Xi_{k}}(t) := \begin{cases} \mathcal{G}_{\hat{\Xi}_{k-1}}(t) & \text{ for } t \in [s, \sigma_{k}] \\ \sup \{ y : D_{BR}( (\mathcal{G}_{\hat{\Xi}_{k-1}}(\sigma_{k}), \sigma_{k}), (y, t)) \leq 1 \} & t> \sigma_{k},\eta_{k} = 1 \\
 \inf \{ y : D_{BR}( (\mathcal{G}_{\hat{\Xi}_{k-1}}(\sigma_{k}), \sigma_{k}), (y, t)) \leq 1 \} & t> \sigma_{k},\eta_{k} = -1 \end{cases}$$
We similarly define functions $G^{n}_{\Xi_{k}} : [\lceil s \rceil_{n}, \infty) \to \mathbb{R}$:
$$G^{n}_{\Xi_{k}}(t) := \begin{cases} G_{\hat{\Xi}_{k-1}}^{n}(t) & \text{ for } t \in [\lceil s \rceil_{n}, \lceil \sigma_{k} \rceil_{n}] \\ \sup \{ y : D_{RW}^{n}( (G^{n}_{\hat{\Xi}_{k-1}}(\lceil \sigma_{k} \rceil_{n}), \lceil \sigma_{k} \rceil_{n}), (y, t)) \leq 1 \} & t \in (\lceil \sigma_{k} \rceil_{n}, \infty) \cap \frac{1}{n}\mathbb{Z},\eta_{k} = 1 \\
\inf \{ y : D_{RW}^{n}( (G^{n}_{\hat{\Xi}_{k-1}}(\lceil \sigma_{k} \rceil_{n}), \lceil \sigma_{k} \rceil_{n}), (y, t)) \leq 1 \} & t \in (\lceil \sigma_{k} \rceil_{n}, \infty) \cap \frac{1}{n}\mathbb{Z}, \eta_{k} = -1 \\ \end{cases}$$
We then extend $G_{\Xi_{k}}^{n}$ to a function defined on $[\lceil s \rceil_{n}, \infty)$ by linearly interpolating between points in $\frac{1}{n} \mathbb{Z} \cap [ \lceil s \rceil_{n}, \infty)$.
\begin{figure}[htp]
    \centering
    \includegraphics[width=15cm]{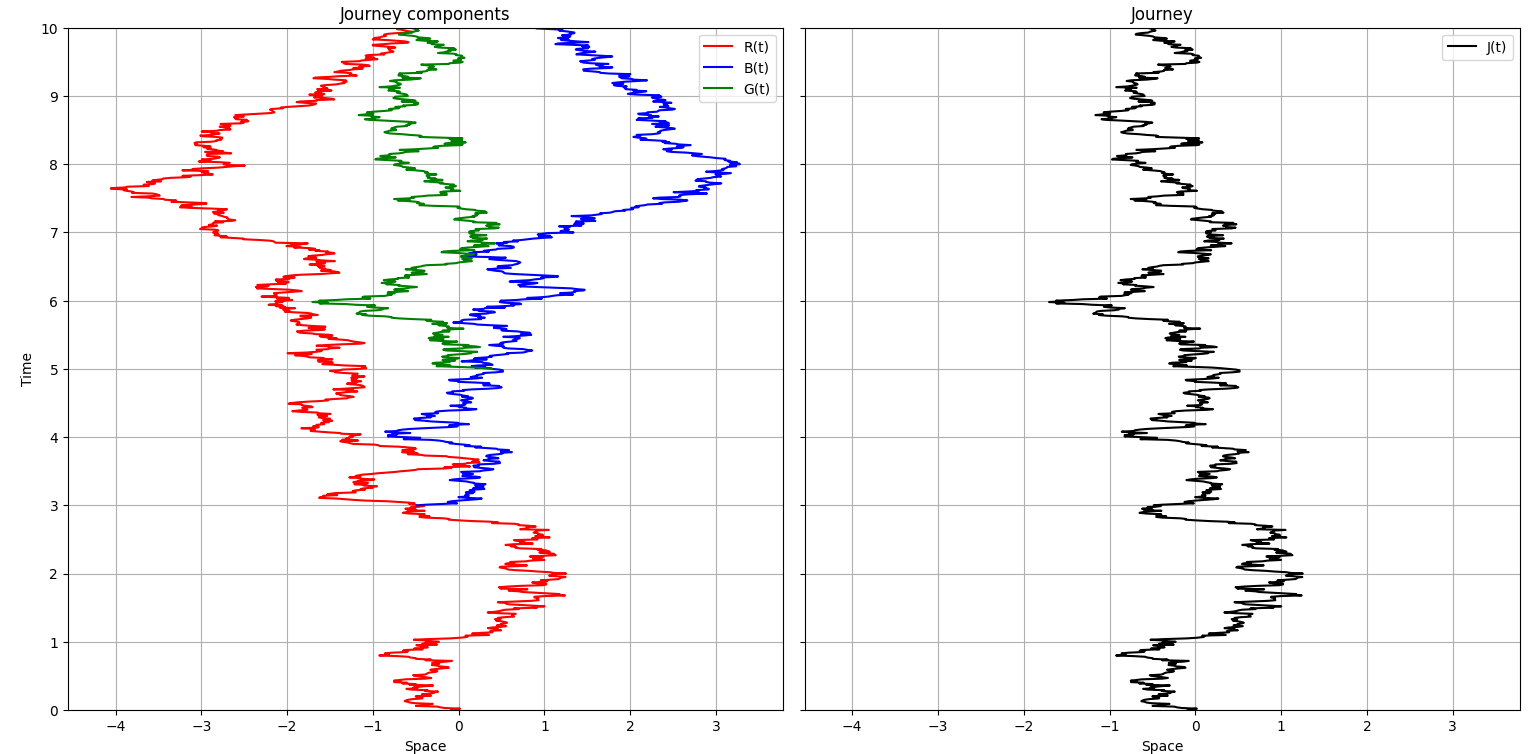}
    \caption{A depiction of a journey $J(t)$, on the Brownian web. The itinerary of $J(t)$ is of length two and takes the values $(x, t) = (0,0)$, $\overline{\sigma} = (3, 5)$ and $\overline{\eta} = (1, -1)$. Since $(x,t)=(0,0)$, we begin by following the path $R(t)$, which is the unique path in $\mathcal{W}$ starting from $(0, 0)$. At time $\sigma=3$, we instead follow $B(t)$. $B(t)$ is the \textbf{right} boundary of all points of distance at most one (according to the distance function $D_{BR}$) from $(R(3), 3)$. The right boundary is chosen as $\eta_{1}=1$. At time $\sigma=5$, we instead follow $G(t)$. $G(t)$ is the \textbf{left} boundary of all points of distance at most one from $(B(5), 5)$. The left boundary is chosen as $\eta_{2}=-1$.}
    \label{fig:journey}
\end{figure} \label{Fig6}
\newpage
\subsection{Preliminaries}
In this subsection, we list some results from other papers used within our arguments.
\newline
\begin{customlemma}{2.1} \label{4.1VV}Let $f^{n}, f$ be lower semicontinuous functions from $\mathbb{R}^{4} \times \overline{\mathbb{R}}$. Suppose that:
\begin{itemize}
\item For any convergent sequence $x^{n} \to x \in \mathbb{R}^{4}$, it holds that:
$$\liminf_{n \to \infty} f^{n}(x^{n}) \geq f(x)$$
\item For any $x \in \mathbb{R}^{4}$, we can find a convergent sequence $x^{n} \to x$ such that:
$$\limsup_{n \to \infty} f^{n}(x^{n}) \leq f(x)$$
\end{itemize}
Then $\mathfrak{e}(f^{n}) \to \mathfrak{e}(f)$ in $(\mathcal{E}_{\star}, d_{\star})$ as $n \to \infty$. 
\begin{proof} See Theorem 4.1 in \cite{VV}.
\end{proof}
\end{customlemma}
\text{ }
\newline

\begin{customlemma}{2.2} \label{L1.2} $\mathcal{W}$ is almost surely a compact subset of $(\Pi, d)$.
\begin{proof} See Proposition 3.2 in \cite{FINR}.
\end{proof}
\end{customlemma}
\text{ }
\newline

\begin{customlemma}{2.3} \label{L1.3} Define the left and right forward-time distance one boundary curves from the point $(x,s)$ as follows:
$$\mathfrak{b}^{R}_{x,s}(t) := \sup \{ y : D_{BR}((x,s), (y,t)) \leq 1 \}$$
$$\mathfrak{b}^{L}_{x,s}(t) := \inf \{ y : D_{BR}((x,s),(y,t)) \leq 1 \}$$
Then $\mathfrak{b}_{x,s}^{R (L)}$ is equal in distribution (as an element of $(\Pi, d)$) to the right (left) Skorokhod reflection (formally defined as $\Lambda_{R}$ $(\Lambda_{L})$ in \hyperref[S5.1]{Section 5.1}) of two independent Brownian motions started at space-time point $(x,s)$.
\begin{proof} This follows from the same argument used in Proposition 3.1 of \cite{VV}. The precise statement of Proposition 3.1 differs from this lemma, however the argument for this lemma is exactly the same.
\end{proof}
\end{customlemma}
\text{ }
\newline

\begin{customlemma}{2.4} \label{L1.4} Let $\{ (\gamma^{n}, s^{n}) \} \subset \mathcal{W}$ be a sequence of paths in $\mathcal{W}$ such that $(\gamma^{n}, s^{n})$ converges to $(\gamma, s) \in \mathcal{W}$ in $(\Pi, d)$ as $n \to \infty$. Let $\sigma^{n}$ be the time at which $\gamma^{n}$ coalesces with $\gamma$. Then
$$\lim_{n \to \infty} \sigma^{n} = s$$
\begin{proof} See Corollary 2.8 in \cite{ESS}.
\end{proof}
\end{customlemma}
\text{ }
\newline
\begin{customlemma}{2.5} \label{L1.5} Almost surely, for all $(x, s) \in \mathbb{R}^{2}$ the following holds:
\newline
\newline
If there exists a path $(\gamma, \hat{s}) \in \mathcal{W}$ with $\hat{s} < s$ and such that $\gamma(s) = x$, then there exist at most two distinct paths in $\mathcal{W}$ originating at the point $(x,s)$. That is to say, if there exists an incoming path to the point $(x,s)$, then there exist at most two distinct outgoing paths from $(x,s)$.
\begin{proof} This follows from the characterization of the special points of the Brownian web, which dates back to \cite{TW}. One proof is given in Theorem 2.11 in \cite{ESS}. 
\end{proof}
\end{customlemma}
\text{ }
\newline

\begin{customlemma}{2.6} \label{L1.6} For $n \in \mathbb{N}$, let $X^{n} := (x_{1}^{n}, x_{2}^{n}, ...)$ and $X := (x_{1}, x_{2}, ...)$ be $\mathbb{R}^{\mathbb{N}}$ valued random variables. Suppose that for $i \in \mathbb{N}$ it is the case that:
$$(x_{1}^{n}, x_{2}^{n}, ..., x_{i}^{n}) \overset{d}{\to} (x_{1}, x_{2}, ..., x_{i})$$
as $n \to \infty$ where convergence takes place in the topology given by the standard product metric. Suppose further that the sequence $\{ X^{n} \}$ is tight. Then:
$$X^{n} \overset{d}{\to} X$$
That is, if $\{ X^{n} \}$ is tight and the finite-dimensional distribution of $X^{n}$ converge to those of $X$, then $X^{n}$ converges weakly to $X$.
\begin{proof} This is proven through multiple theorems in Chapter 1 of \cite{BILL}
\end{proof}
\end{customlemma}
\text{ }
\newline
\begin{customlemma}{2.7} \label{L1.7} $D_{BR}$ is almost surely lower semicontinuous.
\begin{proof} See Theorem 1.4 in \cite{VV}
\end{proof}
\end{customlemma}

\section{Main results}
\subsection{Statement of main results}
Having presented the foundational definitions, we can now describe the main results of the paper:
\begin{customthm}{3.1} \label{T3.1} Suppose that $\mathcal{V}^{n} \overset{d}{\to} \mathcal{W}$ in the space $(\mathcal{H}, d_{\mathcal{H}})$ as $n \to \infty$. Then there exists a coupling of $\mathcal{V}^{n}$ and $\mathcal{W}$ such that the following (almost surely) holds:
$$\forall k \in \mathbb{Z}_{\geq 0}, \text{ all itineraries } \Xi_{k}, \text{ and all } t \in \mathbb{Q} \cap (\sigma_{k}, \infty) \text{ it is the case that }$$ 
$$\lim\limits_{n \to \infty} G^{n}_{\Xi_{k}}(\lceil t \rceil_{n}) = \mathcal{G}_{\Xi_{k}}(t)$$
\end{customthm}
\noindent Theorem \ref{T3.1} is a corollary of Theorem \ref{T5.3.2}. As noted in the introduction, Theorem 1.3.3 in \cite{NRS} provides some simple conditions which can be imposed on $\zeta$ in order to guarantee the distributional convergence of $\mathcal{V}^{n}$ to $\mathcal{W}$. Those results ensure the existence of a broad collection of random variables $\mathcal{V}^{n}$ to which Theorem \ref{T3.1} can be applied. We later leverage Theorem \ref{T3.1} in order to establish convergence of the epigraphs of $D_{RW}^{n}$ to the epigraph of $D_{BR}$:
\begin{customthm}{3.2} \label{T3.2} Under the coupling given in Theorem \ref{T3.1} it is (almost surely) the case that $\mathfrak{e} \left( D_{RW}^{n} \right) \to \mathfrak{e} \left( D_{BR} \right)$ in the space $(\mathcal{E}_{\star}, d_{\star})$ as $n \to \infty$.
\end{customthm}
\noindent Theorem \ref{T3.1} is proven in Section 6 (labeled there as Theorem \ref{T6.2.2} for organizational purposes.) Theorem \ref{T1.1}, which was given in the introduction, follows from Theorem 1.3.3 in \cite{NRS} and Theorem \ref{T3.2}. Although Theorem \ref{T3.1} alone is sufficient to establish the epigraph convergence given in Theorem \ref{T3.2}, the convergence given in Theorem \ref{T3.1} can be strengthened to convergence in $(\Pi, d)$. By showing that the functions $G_{\Xi_{k}}^{n}$ are sufficiently regular, we establish a result showing uniform convergence of journeys along compact sets:
\begin{customthm}{3.3} \label{T3.3} Suppose that $\mathcal{V}^{n} \overset{d}{\to} \mathcal{W}$ in the space $(\mathcal{H}, d_{\mathcal{H}})$ as $n \to \infty$. Then there exists a coupling of $\mathcal{V}^{n}$ and $\mathcal{W}$ such that the following (almost surely) holds:
$$\forall k \in \mathbb{Z}_{\geq 0} \text{ and for all itineraries } \Xi_{k}, \text{ it is the case that } \lim\limits_{n \to \infty} G_{\Xi_{k}}^{n} = \mathcal{G}_{\Xi_{k}} \text{ in } (\Pi, d)$$
\end{customthm}
\noindent Theorem \ref{T3.3} is a corollary of Theorem \ref{T7.3.2}.

\subsection{Outline of the proof of the main results}
\textbf{Section 4:} We begin by developing bounds on the limiting behavior of the distance functions $D_{RW}^{n}$. Loosely speaking, we show that if $\mathcal{V}^{n}$ converges to $\mathcal{W}$, then for fixed points $u, v \in \mathbb{R}^{2}$, the distance $D_{RW}^{n}(u, v)$ is eventually at least as large as the distance $D_{BR}(u, v)$. We use this fact to bound the limiting behavior of the journeys $G_{\Xi_{k}}^{n}$.
\newline
\newline
\textbf{Section 5:} Here we introduce approximating curves $R_{\Xi_{k}}^{n}$. These curves behave similarly to $G_{\Xi_{k}}^{n}$, but have the benefit of being easier to analyze. We establish some convergence results for $R_{\Xi_{k}}^{n}$, and then extend those results to $G_{\Xi_{k}}^{n}$, thereby proving Theorem \ref{T5.3.2}.
\newline
\newline
\textbf{Section 6:} In this section we analyze $\mathcal{W}$ and show that any geodesic in $\mathcal{W}$ can be approximated arbitrarily well by a journey. We use these approximations and Theorem \ref{T5.3.2} to prove Theorem \ref{T6.2.2}.
\newline
\newline
\textbf{Section 7:} Lastly we introduce a modulus of continuity and prove that the functions $R_{\Xi_{k}}^{n}$ are unlikely to vary too wildly. We leverage this continuity result to obtain a similar result for $G_{\Xi_{k}}^{n}$. We then use this continuity result to prove Theorem \ref{T7.3.2}.
\section{Bounds on the limiting behavior of $D_{RW}^{n}$}
Throughout this section, we assume that $\mathcal{V}^{n}$ and $\mathcal{W}$ are defined on some common probability space $\Omega$. We begin with a topological result which provides bounds on the limiting behavior of $D_{RW}^{n}$:
\newtheorem{thm2}{Theorem}
\begin{customthm}{4.1} \label{T4.0.1} On the event that $\mathcal{V}^{n} \to \mathcal{W}$ in $(\mathcal{H}, d_{\mathcal{H}})$, the following almost surely holds:
$$\forall ((x,s),(y,t)) \in \mathbb{R}^{4},$$ $$\lim\limits_{n \to \infty} ((x^{n}, s^{n}), (y^{n}, t^{n})) = ((x,s),(y,t)) \implies$$
$$\liminf_{n \to \infty} D_{RW}^{n}( (x^{n}, s^{n}), (y^{n}, t^{n})) \geq D_{BR}((x,s),(y,t))$$

\begin{proof}
Fix any sequence $\textbf{p}^{n} := ((x^{n}, s^{n}), (y^{n}, t^{n}))$ with $\textbf{p}^{n} \to ((x,s),(y,t))$ and any $\omega \in \Omega$ such that $\mathcal{V}^{n}(\omega) \to \mathcal{W}(\omega)$. Suppose that $\liminf_{n \to \infty} D_{RW}^{n}( \textbf{p}^{n}) = k$. If $k= \infty$ the result trivially holds, so we assume $k \in \mathbb{Z}_{\geq 0}$. Consider some subsequence $\{ n_{j} \} \subseteq \mathbb{N}$ such that $D_{RW}^{n_{j}}(\textbf{p}^{n_{j}}) = k$ for all $j$. 
\newline
\newline
There must exist paths $\pi_{0}^{n_{j}}, ..., \pi_{k}^{n_{j}} \in \mathcal{V}^{n_{j}}$ falling along a geodesic between $(x^{n_{j}}, s^{n_{j}})$ and $(y^{n_{j}}, t^{n_{j}})$. Since $\mathcal{V}^{n} \to \mathcal{W}$ in $(\mathcal{H}, d_{\mathcal{H}})$, for $0 \leq i \leq k$, there exist paths $\gamma^{n_{j}}_{i} \in \mathcal{W}$ such that $d(\gamma^{n_{j}}_{i}, \pi^{n_{j}}_{i}) \to 0$. Note that by Lemma \ref{L1.2}, $\mathcal{W}$ is (a.s.)\ a compact subset of the metric space $(\Pi, d)$. It follows that any sequence of paths in $\mathcal{W}$ has a convergent subsequence. Therefore we can find $\{ n_{j,0} \} \subseteq \{ n_{j} \}$ such that $\gamma^{n_{j, 0}}_{0}$ converges to a limit $\gamma^{\infty}_{0} \in \mathcal{W}$. By repeating the process, we can find a further subsequence $\{n_{j,1} \} \subseteq \{ n_{j,0} \}$ such that $\gamma^{n_{j,1}}_{0}$ converges to $\gamma^{\infty}_{0}$ and $\gamma^{n_{j, 1}}_{1}$ converges to a path, $\gamma^{\infty}_{1} \in \mathcal{W}$. In general, there exists a subsequence $\{ n_{j,k} \} \subseteq \mathbb{N}$ such that it is the case that $\gamma_{i}^{n_{j,k}}$ converge to $\gamma^{\infty}_{i} \in \mathcal{W}$ for all $0 \leq i \leq k$. We consider some such subsequence and the associated limit paths $\{ \gamma^{\infty}_{i}  \}_{0 \leq i \leq k}$.
\newline
\newline
First we claim that $\hat{\gamma}^{\infty}_{k}(t)=y$ and that $\hat{\gamma}^{\infty}_{0}(s)=x$. Let $j$ be large enough so that $|s^{n_{j, k}} - s| < 1$ and that $|t^{n_{j, k}} - t| < 1$. Let $s_{i}^{n_{j,k}}$ denote the starting time of path $\pi_{i}^{n_{j,k}}$. Let $w_{i}^{n_{j, k}}$ be the time at which path $i$ jumps to path $i+1$. Then for all $i$ and all large enough $j$, we have: $w_{i}^{n_{j, k}}, s_{i}^{n_{j, k}} \in [s-1, t+1]$. Therefore by Lemma \ref{L9.0.1}, it is the case that $\gamma_{k}^{\infty}, \gamma_{0}^{\infty} \in \Pi^{\star}$, so that $\hat{\gamma}_{k}^{\infty}, \hat{\gamma}_{0}^{\infty}$ are continuous with respect to the standard topology on $\mathbb{R}$. By Lemma \ref{L9.0.2} it follows that $\hat{\pi}_{k}^{n_{j, k}}$ converges uniformly to $\hat{\gamma}^{\infty}_{k}$ on $[t-1, t+1]$. Using uniform convergence and the continuity of $\hat{\gamma}_{k}^{\infty}$, we have:
$$y = \lim\limits_{j \to \infty} y^{n_{j, k}} = \lim\limits_{j \to \infty} \hat{\pi}_{k}^{n_{j, k}}(t^{n_{j, k}}) = \lim\limits_{j \to \infty} \hat{\gamma}_{k}^{\infty}(t^{n_{j, k}}) = \hat{\gamma}_{k}^{\infty}(t)$$
Similarly, since $\hat{\pi}^{n_{j, k}}_{0}$ converges uniformly to $\hat{\gamma}_{0}^{\infty}$ on $[s-1, s+1]$, we have:
$$x = \lim\limits_{j \to \infty} x^{n_{j, k}} = \lim\limits_{j \to \infty} \hat{\pi}_{0}^{n_{j, k}}(s^{n_{j, k}}) = \lim\limits_{j \to \infty} \gamma_{0}^{\infty} (s^{n_{j, k}}) = \hat{\gamma}_{0}^{\infty}(s)$$
Next we claim that, if the paths $\gamma^{\infty}_{i}$ start at times $s_{i}$ then we have the following ordering: $$s = s_{0} \leq s_{1} \leq ... \leq s_{k} \leq t$$
This follows because convergence in $(\Pi, d)$ implies convergence of start times, hence:
$$s_{i}-s_{i-1} = \lim\limits_{j \to \infty} \left( s_{i}^{n_{j, k}}-s_{i-1}^{n_{j, k}} \right) \geq 0$$
$$s_{0} = \lim\limits_{j \to \infty} s_{0}^{n_{j, k}} = s$$
$$s_{k} = \lim\limits_{j \to \infty} s_{k}^{n_{j, k}} \leq \lim\limits_{j \to \infty} t^{n_{j, k}} = t$$
\newline
Next we claim that for $0 < i \leq k$, it is the case that $\gamma^{\infty}_{i}(s_{i}) = \gamma^{\infty}_{i-1}(s_{i})$. Let 
$$c:= \max \{|y| : (y,w) \in \mathcal{K} \}$$
$$h := \max \{|w|: (y, w) \in \mathcal{K} \}$$
By definition of the jumps in our model $D_{RW}^{n}$, the following hold:
$$|\pi_{i-1}(w_{i}^{n_{j, k}}) - \pi_{i}(s_{i}^{n_{j, k}})| \leq \frac{c}{\sqrt{n_{j, k}}}$$
$$|w_{i}^{n_{j, k}} - s_{i}^{n_{j, k}}| \leq \frac{h}{n_{j, k}}$$
By the triangle inequality we have:
\begin{equation*}
\begin{split}
|\gamma_{i}^{\infty}(s_{i})-\gamma_{i-1}^{\infty}(s_{i})| = |\hat{\gamma}_{i}^{\infty}(s_{i}) -\hat{\gamma}_{i-1}^{\infty}(s_{i})| &\leq |\hat{\gamma}_{i}^{\infty}(s_{i})-\hat{\gamma}_{i}^{\infty}(s_{i}^{n_{j, k}})| \\
&+ |\hat{\gamma}_{i}^{\infty}(s_{i}^{n_{j, k}}) - \hat{\pi}_{i}^{n_{j, k}}(s_{i}^{n_{j, k}})| \\
&+ |\hat{\pi}^{n_{j, k}}_{i}(s_{i}^{n_{j, k}}) - \hat{\pi}^{n_{j, k}}_{i-1}(w_{i}^{n_{j, k}})| \\
&+ |\hat{\pi}_{i-1}^{n_{j, k}}(w_{i}^{n_{j, k}}) - \hat{\gamma}_{i-1}^{\infty}(w_{i}^{n_{j, k}})| \\
&+ |\hat{\gamma}_{i-1}^{\infty}(w_{i}^{n_{j, k}})-\hat{\gamma}_{i-1}^{\infty}(s_{i})|
\end{split}
\end{equation*}
Since $\lim_{j \to \infty} s_{i}^{n_{j, k}} = s_{i}$ and $\hat{\gamma}_{i}^{\infty}$ are continuous, for any $\epsilon>0$ and all large enough $j$ we have:
$$|\hat{\gamma}_{i}^{\infty}(s_{i})-\hat{\gamma}_{i}^{\infty}(s_{i}^{n_{j, k}})| < \epsilon$$
Since $|w_{i}^{n_{j, k}} - s_{i}^{n_{j, k}}| \leq \frac{h}{n_{j, k}}$, we also have $\lim_{j \to \infty} w_{i}^{n_{j, k}} = s_{i}$. For large enough $j$ we then also have:
$$|\hat{\gamma}_{i-1}^{\infty}(w_{i}^{n_{j, k}})-\hat{\gamma}_{i-1}^{\infty}(s_{i})| < \epsilon$$
Since $\hat{\pi}_{i}^{n_{j, k}}$ converges uniformly to $\hat{\gamma}_{i}^{\infty}$, for large enough $j$ we also have:
$$|\hat{\pi}^{n_{j, k}}_{i}(s_{i}^{n_{j, k}}) - \hat{\gamma}_{i}^{\infty}(s_{i}^{n_{j, k}})| \leq \sup\limits_{r \in [s-1, t+1]} |\hat{\pi}^{n_{j, k}}_{i}(r) - \hat{\gamma}^{\infty}_{i}(r) | < \epsilon$$
The same argument shows for large enough $j$ we have:
$$|\hat{\pi}_{i-1}^{n_{j, k}}(w_{i}^{n_{j, k}}) - \hat{\gamma}_{i-1}^{\infty}(w_{i}^{n_{j, k}})| < \epsilon$$
Taking $j$ large enough such that $\frac{c}{\sqrt{n_{j, k}}} < \epsilon$, we also have:
$$|\hat{\pi}^{n_{j, k}}_{i-1}(w_{i}^{n_{j, k}}) - \hat{\pi}^{n_{j, k}}_{i}(s_{i}^{n_{j, k}})| \leq \frac{c}{\sqrt{n_{j, k}}} < \epsilon$$
Altogether, this shows that for any $\epsilon > 0$, it is the case that $|\gamma_{i}^{\infty}(s_{i})-\gamma_{i-1}^{\infty}(s_{i})| < 5 \epsilon$. Hence $|\gamma_{i}^{\infty}(s_{i})-\gamma_{i-1}^{\infty}(s_{i})|=0$ as desired.
\newline
\newline
We have shown that if $\mathcal{V}^{n}(\omega) \to \mathcal{W}(\omega)$ then there (a.s.) exist a sequence of continuous paths, $\gamma_{0}^{\infty}, ..., \gamma_{k}^{\infty} \in \mathcal{W}$ satisfying:
\begin{enumerate}[label=(\roman*)]
        \item $\gamma_{0}^{\infty}(s) = x$
        \item $\gamma_{k}^{\infty}(t) = y$
        \item $\gamma_{i}^{\infty}$ starts at time $s_{i}$
        \item $s = s_{0} \leq ... \leq s_{k} \leq t$
	\item For $1 \leq i \leq k$, $\gamma_{i}^{\infty}(s_{i}) = \gamma_{i-1}^{\infty}(s_{i})$
    \end{enumerate}
It follows that $D_{BR}((x,s),(y,t)) \leq k$.
\end{proof}
\end{customthm}
\text{ }
\newline
We now leverage Theorem \ref{T4.0.1} to give bounds on the pointwise limits of our journeys.
\newline
\begin{customthm}{4.2} \label{T4.0.2} Let $k \in \mathbb{Z}_{\geq 0}$ and let $\Xi_{k} = (x,s, \overline{\sigma}, \overline{\eta})$ be a length $k$ itinerary. Assume the following hold:
\begin{itemize}
\item $\mathcal{V}^{n} \to \mathcal{W}$ in $(\mathcal{H}, d_{\mathcal{H}})$
\item $G_{\Xi_{k}}^{n}(\lceil \sigma_{k} \rceil_{n}) \to \mathcal{G}_{\Xi_{k}}(\sigma_{k}) \text{ if } k \geq 1$
\end{itemize}
Then the following (a.s.) holds: 
\newline
\newline
For any $w > \sigma_{k}$ (with $w>s$ for $k=0$) and any sequence $\{ w^{n} \}$ such that $\lim_{n \to \infty} w^{n}=w$:
\begin{itemize}
\item If $\eta_{k} = 1$ or $k=0$ then $\limsup\limits_{n \to \infty} G_{\Xi_{k}}^{n}(w^{n}) \leq \mathcal{G}_{\Xi_{k}}(w)$
\item If $\eta_{k} = -1$ or $k=0$ then $\liminf\limits_{n \to \infty} G_{\Xi_{k}}^{n}(w^{n}) \geq \mathcal{G}_{\Xi_{k}}(w)$
\end{itemize}

\begin{proof} We prove only the cases where $\eta_{k}=1$ or $k=0$, the case where $\eta_{k}=-1$ is symmetric to the case where $\eta_{k}=1$. First we deal with the case $k=0$. Let $\{n_{j} \} \subseteq \mathbb{N}$. By the assumed convergence of $\mathcal{V}^{n} \to \mathcal{W}$, there must exist $\gamma^{n_{j}} \in \mathcal{W}$ such that:
$$\lim_{j \to \infty} d(\gamma^{n_{j}}, Y_{\lceil x \rceil_{\sqrt{n_{j}}}, \lceil t \rceil_{n_{j}}}^{n_{j}}) = 0$$
By the compactness of $\mathcal{W}$, there must exist a subsequence $\{ n_{j, 0} \} \subseteq \mathbb{N}$ and limiting curve $\gamma^{\infty} \in \mathcal{W}$ such that:
$$\lim_{n \to \infty} d(\gamma^{n_{j, 0}}, \gamma^{\infty}) = 0$$
Since $\gamma^{\infty}$ must start from the space-time point $\lim_{j \to \infty} (\lceil x \rceil_{\sqrt{n_{j, 0}}}, \lceil t \rceil_{n_{j, 0}})=(x,t) \in \mathbb{Q}^{2}$, almost surely there is a unique such $\gamma^{\infty}$, and it is equal to $\mathcal{G}_{\Xi_{0}}$ . It follows that any subsequence $\{ n_{j} \} \subseteq \mathbb{N}$ must have a further subsequence $\{ n_{j, 0} \} \subseteq \{ n_{j} \}$ along which $Y^{n_{j, 0}}_{\lceil x \rceil_{\sqrt{n_{j, 0}}}, \lceil t \rceil_{n_{j, 0}}}$ converges to $\gamma^{\infty}$. It follows that $Y^{n}_{\lceil x \rceil_{\sqrt{n}}, \lceil t \rceil_{n}}$ converges to $\gamma^{\infty}$. Then by Lemma \ref{L9.0.1}, \ref{L9.0.2}, we have:
\begin{equation*}
\begin{split}
\lim_{n \to \infty} G^{n}_{\Xi_{0}}(w^{n}) &= \lim_{n \to \infty} Y^{n}_{\lceil x \rceil_{\sqrt{n}}, \lceil t \rceil_{n}}(w^{n}) \\
&= \lim_{n \to \infty} \mathcal{G}_{\Xi_{0}}(w^{n}) \\
&= \mathcal{G}_{\Xi_{0}}(w)
\end{split}
\end{equation*}
We now move onto the case where $k \geq 1$. Let:
$$w_{\star}^{n} := \argmax_{ t \in \{ \lceil w^{n} \rceil_{n}, \lfloor w^{n} \rfloor_{n} \} } G_{\Xi_{k}}^{n}(t)$$
Suppose that $\eta_{k}=1$ and $\limsup_{n \to \infty} G_{\Xi_{k}}^{n}(w^{n}_{\star}) = y \in \overline{\mathbb{R}}$. If $y = - \infty$ the inequality trivially holds as $G_{\Xi_{k}}^{n}$ is linearly interpolated between consecutive time points. Suppose instead $y \in \mathbb{R}$. Then there must exist a subsequence $\{ n_{j} \} \subseteq \mathbb{N}$ such that:
$$\lim_{j \to \infty} G_{\Xi_{k}}^{n_{j}}(w_{\star}^{n_{j}}) = y$$
Then define:
$$y^{n} := \begin{cases} G_{\Xi_{k}}^{n_{j}}(w_{\star}^{n_{j}}) & \text{ for } n \in \{n_{j} \} \\ y & \text{ else } \end{cases}$$
We would have:
$$\liminf\limits_{n \to \infty} D_{RW}^{n}( (G_{\Xi_{k}}^{n}( \lceil \sigma_{k} \rceil_{n}), \lceil \sigma_{k} \rceil_{n}), (y^{n}, w_{\star}^{n})) \leq 1$$
By assumption, $G_{\Xi_{k}}^{n}( \lceil \sigma_{k} \rceil_{n}) \to \mathcal{G}_{\Xi_{k}}(\sigma_{k})$, hence:
$$\lim\limits_{n \to \infty} ((G_{\Xi_{k}}^{n}( \lceil \sigma_{k} \rceil_{n}), \lceil \sigma_{k} \rceil_{n}), (y^{n}, w_{\star}^{n})) = ((\mathcal{G}_{\Xi_{k}}(\sigma_{k}), \sigma_{k}),(y,w))$$
Then, by Theorem \ref{T4.0.1} we would have:
$$D_{BR}((\mathcal{G}_{\Xi_{k}}(\sigma_{k}), \sigma_{k}),(y,w)) \leq 1$$
Hence showing $\mathcal{G}_{\Xi_{k}}(w) \geq y$ as desired.
\newline
\newline
Lastly, suppose $y = \infty$. Take some subsequence $\{n_{j}\}$ along which $\lim_{j \to \infty} G_{\Xi_{k}}^{n_{j}}(w_{\star}^{n_{j}}) = \infty$. Let $\pi_{1}^{n_{j}} \in \mathcal{V}^{n_{j}}$ be the path starting at the point $(G_{\Xi_{k}}^{n_{j}}( \lceil \sigma_{k} \rceil_{n_{j}}), \lceil \sigma_{k} \rceil_{n_{j}})$ and let $\pi_{2}^{n_{j}}$ be the second path along a geodesic between $(G_{\Xi_{k}}^{n_{j}}( \lceil \sigma_{k} \rceil_{n_{j}}), \lceil \sigma_{k} \rceil_{n_{j}})$ and $(G_{\Xi_{k}}^{n_{j}}(w_{\star}^{n_{j}}), w_{\star}^{n_{j}})$. As argued in \ref{T4.0.1}, since $\mathcal{V}^{n} \to \mathcal{W}$, there must exist paths $\gamma_{1}^{\infty}, \gamma_{2}^{\infty} \in \mathcal{W}$ and a subsequence $\{ n_{j, 2} \} \subseteq \{ n_{j} \}$ such that $\pi_{1}^{n_{j,2}}$ converges to $\gamma_{1}^{\infty}$ and $\pi_{2}^{n_{j,2}}$ converges to $\gamma_{2}^{\infty}$. By \ref{L9.0.1}, we would have $\gamma_{1}^{\infty}, \gamma_{2}^{\infty} \in \Pi^{\star}$. Since $\gamma_{2}^{\infty} \in \Pi^{\star}$ we would have: $$\sup\limits_{t \in [w-1, w+1]} \hat{\gamma}_{2}^{\infty}(t) < \infty$$
Due to the uniform convergence given by Lemma \ref{L9.0.2}, this would contradict the fact that $\lim_{j \to \infty} \pi_{2}^{n_{j,2}}(w^{n_{j,2}}) = \infty$. Hence $y \neq \infty$. 
\end{proof}
\end{customthm}
\section{Pointwise convergence of journeys}
\subsection{Definition of approximation curves} \label{S5.1}
\begin{figure}[htp] \label{Fig7}
    \centering
    \includegraphics[width=15cm]{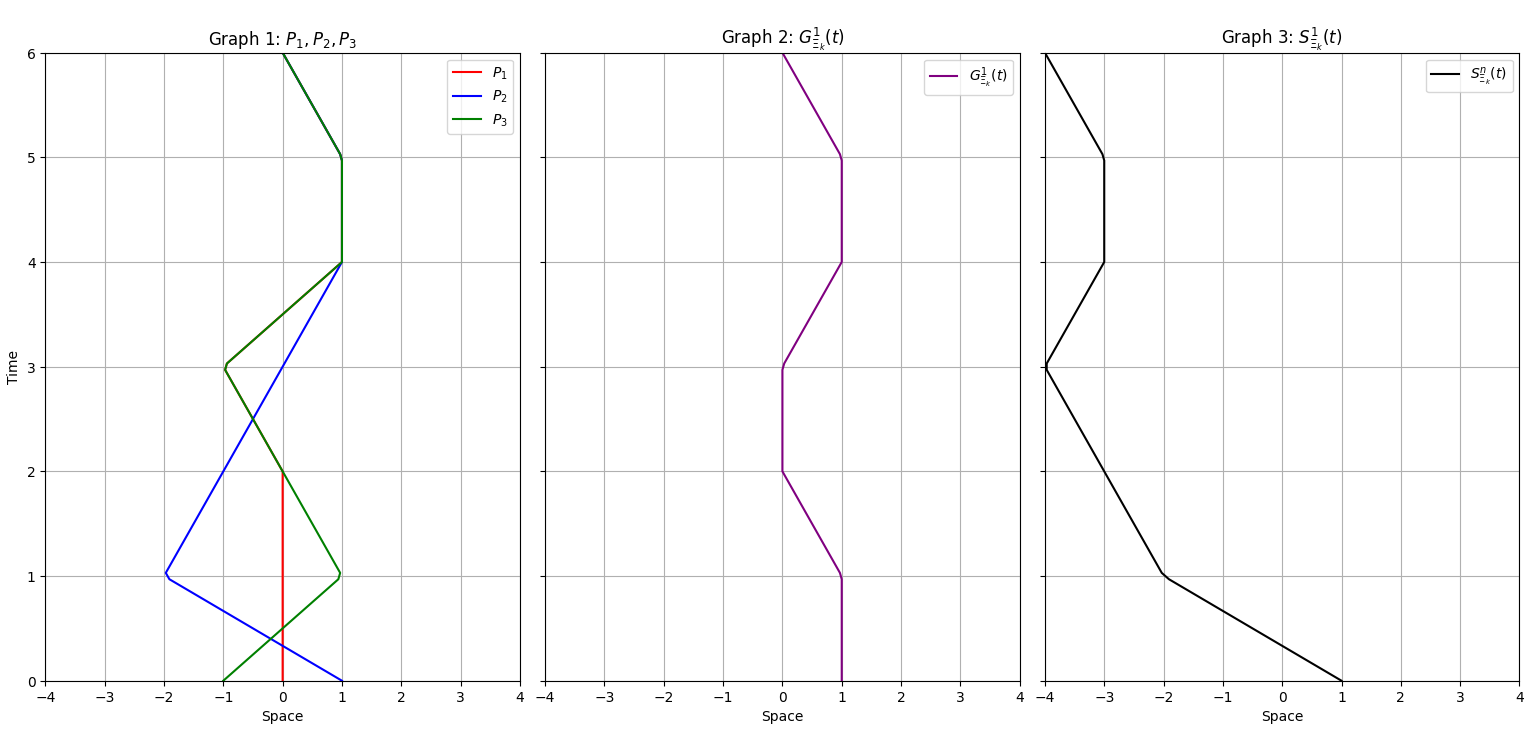}
    \caption{Here $\eta_{k} = 1$. $P_{1}, P_{2}, P_{3} \in \mathcal{V}^{1}$ represent the three right-most paths of distance one from $(G_{\Xi_{k}}^{1}( \lceil \sigma_{k} \rceil ), \lceil \sigma_{k} \rceil)$ where $\sigma_{k} <0$ is the $k$'th jump time, occurring outside the time window displayed in the graph. At integer times $t$, we see that $G_{\Xi_{k}}^{1}(t) = \max \{ P_{1}(t), P_{2}(t), P_{3}(t) \}$. The evolution of $S_{\Xi_{k}}^{1}$ is distributed as a random walk with increment $\zeta$ -- it has an increment which matches the increment of the right-most path, however it does not necessarily evolve in the same way as $G_{\Xi_{k}}^{n}$. Once the three paths coalesce, the evolution of $S_{\Xi_{k}}^{n}$ matches that of $G_{\Xi_{k}}^{n}$ until another path of distance one from $(G_{\Xi_{k}}^{1}(\lceil \sigma_{k} \rceil), \lceil \sigma_{k} \rceil)$ crosses to the right of the paths $P_{i}$.}
    \label{fig:srep}
\end{figure}

In this section we define auxiliary functions $R_{\Xi_{k}}^{n}$ and $S_{\Xi_{k}}^{n}$ which can be used to analyze the behavior of $G_{\Xi_{k}}^{n}$. 
\newline
\newline
Let $k \geq 1$. For $m \in \mathbb{Z}_{\geq 0}$, let $\alpha_{k, m}^{n} := \lceil \sigma_{k} \rceil_{n} + \frac{m}{n}$ and let $\beta_{k,m}^{n} := G_{\Xi_{k}}^{n}(\alpha_{k, m}^{n})$. Define a function $S_{\Xi_{k}}^{n} : [ \lceil \sigma_{k} \rceil_{n}, \infty) \to \mathbb{R}$ satisfying the following relations:
$$S_{\Xi_{k}}^{n}( \lceil \sigma_{k} \rceil_{n}) = G_{\Xi_{k}}^{n}( \lceil \sigma_{k} \rceil_{n})$$
$$S_{\Xi_{k}}^{n}(\alpha_{k, m+1}^{n}) - S_{\Xi_{k}}^{n}(\alpha^{n}_{k, m}) = Y^{n}_{\beta_{k, m}^{n}, \alpha_{k, m}^{n}}(\alpha^{n}_{k, m+1})-Y^{n}_{\beta_{k, m}^{n}, \alpha_{k, m}^{n}}(\alpha^{n}_{k, m})$$
With $S_{\Xi_{k}}^{n}$ linear interpolated between points on $[ \lceil \sigma_{k} \rceil_{n}, \infty) \cap \frac{1}{n}\mathbb{Z}$. Intuitively, $S_{\Xi_{k}}^{n}$ is a random walk whose increment matches the increment of the path located at $(\beta_{k, m}^{n}, \alpha_{k, m}^{n})$. Note that this does not necessarily describe the evolution of $G_{\Xi_{k}}^{n}$ itself (see \hyperref[Fig7]{Figure 7}.)

\begin{figure}[htp] \label{Fig8}
    \centering
    \includegraphics[width=15cm]{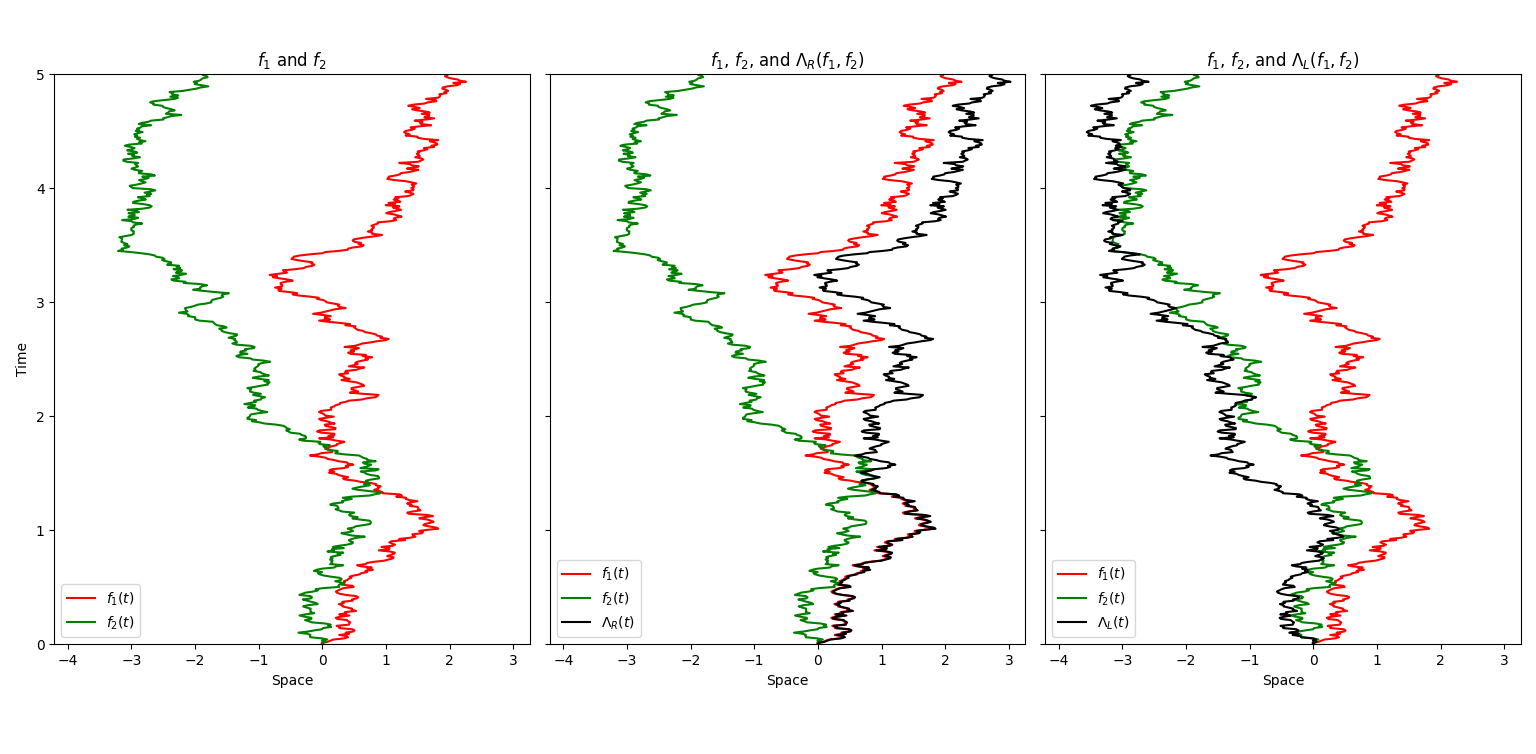}
    \caption{Here $f_{1}(t)$ and $f_{2}(t)$ are chosen to be two independent Brownian motions started at $(0, 0)$. $\Lambda_{R}(f_{1}, f_{2})$ evolves in the same way as $f_{1}$, unless doing so would cause it to cross to the left side of $f_{2}(t)$. As a result, in the case that the $f_{1}, f_{2}$ share a starting point, $\Lambda_{R}(f_{1}, f_{2})$ always lies on the right side of both $f_{1}$ and $f_{2}$. Similarly $\Lambda_{L}(f_{1}, f_{2})$ evolves in the same way as $f_{1}$ unless doing so would result in it crossing to the right of $f_{2}(t)$. In the case that $f_{1}, f_{2}$ share a starting point, $\Lambda_{L}(f_{1}, f_{2})$ always lies to the left of $f_{1}$ and $f_{2}$.}
    \label{fig:refmap}
\end{figure}
\text{ }
\newline
\noindent In order to define the approximation curves, $R_{\Xi_{k}}^{n}$, we introduce two reflection maps, first the map $\Lambda_{R} : (\Pi^{\star})^{2} \to \Pi^{\star}$ defined as follows:
$$\Lambda_{R}((f_{1}, t_{1}), (f_{2}, t_{2})) := (h, t_{3})$$
$$t_{3} := t_{1} \vee t_{2}$$
$$h(t) := \hat{f}_{1}(t) - \inf\limits_{t_{3} \leq s \leq t} \left[ \hat{f}_{1}(s) - \hat{f}_{2}(s) \right]$$
Then the map  $\Lambda_{L} : (\Pi^{\star})^{2} \to \Pi^{\star}$, defined as:
$$\Lambda_{L}((f_{1}, t_{1}), (f_{2}, t_{2})) := (h, t_{3})$$
$$t_{3} := t_{1} \vee t_{2}$$
$$h(t) := \hat{f}_{1}(t) - \sup\limits_{t_{3} \leq s \leq t} \left[ \hat{f}_{1}(s) - \hat{f}_{2}(s) \right]$$
$\Lambda_{R}$ and $\Lambda_{L}$ describe the right and left Skorokhod reflections of the path $f_{1}$ off of the path $f_{2}$. (\hyperref[Fig8]{Figure 8})
\newline
\newline
\noindent We now use the maps $\Lambda_{R}, \Lambda_{L}$ to define our approximation curves $R_{\Xi_{k}}^{n}$. Let $\Xi_{k} := (x, s, \overline{\sigma}, \overline{\eta})$. Let $Y^{n}_{\beta_{k, 0}^{n}, \lceil \sigma_{k} \rceil_{n}} \in \mathcal{V}^{n}$ be the path started at point $(\beta_{k, 0}^{n}, \lceil \sigma_{k} \rceil_{n})$. Then for $k \geq 1$, define the approximation curves:
$$R_{\Xi_{k}}^{n} := \begin{cases} \Lambda_{R}\left( S_{\Xi_{k}}^{n}, Y^{n}_{\beta_{k, 0}^{n}, \lceil \sigma_{k} \rceil_{n}} \right) & \text{ if } \eta_{k} = 1 \\ \Lambda_{L}\left( S_{\Xi_{k}}^{n}, Y^{n}_{\beta_{k, 0}^{n}, \lceil \sigma_{k} \rceil_{n}} \right) & \text{ if } \eta_{k} = -1 \end{cases}$$

\noindent  In words, $R_{\Xi_{k}}^{n}$ is the right (left) Skorokhod reflection of $S_{\Xi_{k}}^{n}$ off of $Y_{\beta_{k, 0}^{n}, \lceil \sigma_{k} \rceil_{n}}^{n}$.
\newline
\newline
The intuition behind choosing this definitions or our approximation curves can be seen through the results of Lemma \ref{L1.3}. The forward-time boundary curves of a point in the Brownian web are distributed as the Skorokhod reflection of two Brownian motions, so we define our discrete approximations using the Skorokhod reflection of a discrete random walk off of another discrete random walk.

\subsection{Properties of approximation curves}
We now prove a number of results relating to the auxiliary functions $S_{\Xi_{k}}^{n}$ and $R_{\Xi_{k}}^{n}$. Lemma \ref{L5.2.1} and Lemma \ref{L5.2.2} show that the approximation curves $R_{\Xi_{k}}^{n}$ are in some way dominated by the journeys $G_{\Xi_{k}}^{n}$; that is to say, $G_{\Xi_{k}}^{n}(t)$ always falls to the right (left in the case that $\eta_{k}=-1$) of $R_{\Xi_{k}}^{n}$. Lemma \ref{L5.2.3} and Lemma \ref{L5.2.4} show that the approximation curves $R_{\Xi_{k}}^{n}$ converge in distribution to a reflected pair of Brownian motions.

\begin{customlemma}{5.1} \label{L5.2.1}
Let $k \geq 1$ and let $\lceil \sigma_{k} \rceil_{n} \leq w \leq t$. Then:
$$S_{\Xi_{k}}^{n}(t) - S_{\Xi_{k}}^{n}(w) \leq G_{\Xi_{k}}^{n}(t) - G_{\Xi_{k}}^{n}(w) \text{ if } \eta_{k} = 1$$
$$S_{\Xi_{k}}^{n}(t) - S_{\Xi_{k}}^{n}(w) \geq G_{\Xi_{k}}^{n}(t) - G_{\Xi_{k}}^{n}(w) \text{ if } \eta_{k} = -1$$
\begin{proof} We prove only the first result, as the argument for the second is symmetric. Let $\eta_{k}=1$. Let $m \in \mathbb{Z}_{\geq 0}$, first we claim the result holds for consecutive lattice points, $\alpha_{k, m}^{n}, \alpha_{k,m+1}^{n}$. Note that:
$$Y^{n}_{\beta_{k, m}^{n}, \alpha_{k, m}^{n}}(\alpha^{n}_{k, m}) = G_{\Xi_{k}}^{n}(\alpha_{k, m}^{n}) \implies$$
$$D_{RW}^{n}((G_{\Xi_{k}}^{n}( \lceil \sigma_{k} \rceil_{n}), \lceil \sigma_{k} \rceil_{n}), (Y^{n}_{\beta_{k, m}^{n}, \alpha_{k, m}^{n}}(\alpha^{n}_{k, m}), \alpha^{n}_{k, m})) \leq 1 \implies$$ 
$$D_{RW}^{n}((G_{\Xi_{k}}^{n}( \lceil \sigma_{k} \rceil_{n}), \lceil \sigma_{k} \rceil_{n}), (Y^{n}_{\beta_{k, m}^{n}, \alpha_{k, m}^{n}}(\alpha^{n}_{k, m+1}), \alpha^{n}_{k, m+1})) \leq 1 \implies$$ 
$$G_{\Xi_{k}}^{n}(\alpha_{k, m+1}^{n}) \geq Y^{n}_{\beta_{k, m}^{n}, \alpha_{k, m}^{n}}(\alpha^{n}_{k, m+1})$$
Hence:
\begin{equation*}
\begin{split}
S_{\Xi_{k}}^{n}(\alpha_{k, m+1}^{n}) - S_{\Xi_{k}}^{n}(\alpha^{n}_{k, m}) &= Y^{n}_{\beta_{k, m}^{n}, \alpha_{k, m}^{n}}(\alpha^{n}_{k, m+1})-Y^{n}_{\beta_{k, m}^{n}, \alpha_{k, m}^{n}}(\alpha^{n}_{k, m}) \\
&\leq G_{\Xi_{k}}^{n}(\alpha_{k, m+1}^{n}) -Y^{n}_{\beta_{k, m}^{n}, \alpha_{k, m}^{n}}(\alpha^{n}_{k, m}) \\
&= G_{\Xi_{k}}^{n}(\alpha_{k, m+1}^{n}) -G_{\Xi_{k}}^{n}(\alpha_{k, m}^{n})
\end{split}
\end{equation*}
Thus the result holds for consecutive lattice points, hence holds for $w,t$ with $\lceil \sigma_{k} \rceil_{n} \leq w \leq t$ with $w, t \in \frac{1}{n} \mathbb{Z}$. As both $S_{\Xi_{k}}^{n}$ and $G_{\Xi_{k}}^{n}$ are defined to be the linear interpolation of such points, it then holds for all $w,t$ with $\lceil \sigma_{k} \rceil_{n} \leq w \leq t$.
\end{proof}
\end{customlemma}
\text{ }
\newline
\begin{customlemma}{5.2} \label{L5.2.2}
For any $k \geq 1$ and $t \geq \lceil \sigma_{k} \rceil_{n}$, it is the case that:
$$R_{\Xi_{k}}^{n}(t) \leq G_{\Xi_{k}}^{n}(t) \text{ if } \eta_{k} = 1$$
$$R_{\Xi_{k}}^{n}(t) \geq G_{\Xi_{k}}^{n}(t) \text{ if } \eta_{k} = -1$$
\begin{proof} 
As above, we prove only the first result, the second follows from a symmetric argument. Suppose $\eta_{k} = 1$. Note that for $w \in [ \lceil \sigma_{k} \rceil_{n}, \infty) \cap \frac{1}{n} \mathbb{Z}$ we have:
$$D_{RW}^{n}((G_{\Xi_{k}}^{n}( \lceil \sigma_{k} \rceil_{n}), \lceil \sigma_{k} \rceil_{n}), (Y^{n}_{\beta_{k, 0}^{n}, \lceil \sigma_{k} \rceil_{n}}(w), w)) = 0 \leq 1$$
Hence, by the fact that $G_{\Xi_{k}}^{n}, Y^{n}_{\beta_{k, 0}^{n}, \lceil \sigma_{k} \rceil_{n}}$ are linear in between such points, we have $G_{\Xi_{k}}^{n}(w) \geq Y^{n}_{\beta_{k, 0}^{n}, \lceil \sigma_{k} \rceil_{n}}(w)$ for all $w \in [ \lceil \sigma_{k} \rceil_{n}, \infty)$.
Define:
$$\tau := \sup\limits_{w \in [\lceil \sigma_{k} \rceil_{n}, t]} \{ w : R_{\Xi_{k}}^{n}(w) =Y^{n}_{\beta_{k, 0}^{n}, \lceil \sigma_{k} \rceil_{n}}(w) \}$$
By definition of $\Lambda_{R}$, $R_{\Xi_{k}}^{n}$ evolves as $S_{\Xi_{k}}^{n}$ along $[\tau, t]$. Then by Lemma \ref{L5.2.1}, we have:
\begin{equation*}
\begin{split}
R_{\Xi_{k}}^{n}(t) &= \left[ R_{\Xi_{k}}^{n}(t)-R_{\Xi_{k}}^{n}(\tau) \right] +R_{\Xi_{k}}^{n}(\tau) \\
&= \left[ S_{\Xi_{k}}^{n}(t) - S_{\Xi_{k}}^{n}(\tau) \right] +Y^{n}_{\beta_{k, 0}^{n}, \lceil \sigma_{k} \rceil_{n}}(\tau) \\
& \leq \left[ G_{\Xi_{k}}^{n}(t) - G_{\Xi_{k}}^{n}(\tau) \right] + G_{\Xi_{k}}^{n}(\tau) \\
&= G_{\Xi_{k}}^{n}(t)
\end{split}
\end{equation*}
\end{proof}
\end{customlemma}
\begin{customlemma}{5.3} \label{L5.2.3}
 Let $k \geq 1$ and $t > \alpha_{k,1}^{n}$. Then $S_{\Xi_{k}}^{n}(t) - S_{\Xi_{k}}^{n}(\alpha_{k,1}^{n})$ is independent of $Y^{n}_{\beta_{k, 0}^{n}, \lceil \sigma_{k} \rceil_{n}}(t) -Y^{n}_{\beta_{k, 0}^{n}, \lceil \sigma_{k} \rceil_{n}}(\alpha_{k, 1}^{n})$.
\begin{proof} Without loss of generality, we let $\eta_{k}=1$. The evolution of $S_{\Xi_{k}}^{n}$ is determined by increment variables $\zeta$ associated with the vertices $(\beta^{n}_{k, m}, \alpha^{n}_{k, m})$. Since $(x, 0) \in \mathcal{K}$ for some $x>0$, it follows that for $m \geq 1$:
$$\beta^{n}_{k, m} = G_{\Xi_{k}}^{n}(\alpha_{k, m}^{n}) \geq Y^{n}_{\beta_{k, 0}^{n},  \lceil \sigma_{k} \rceil_{n}}(\alpha_{k, m}^{n}) + x > Y^{n}_{\beta_{k, 0}^{n}, \lceil \sigma_{k} \rceil_{n}}(\alpha_{k, m}^{n})$$
Hence for all $m \geq 1$, it is the case that $\beta^{n}_{k, m} > Y^{n}_{\beta_{k, 0}^{n}, \lceil \sigma_{k} \rceil_{n}}(\alpha_{k, m}^{n})$. So the set of increment variables defining the evolution of $S_{\Xi_{k}}^{n}$ along $[\alpha_{k, 1}^{n}, \infty)$ is disjoint from the set of increment variables defining the evolution of $Y^{n}_{\beta^{n}_{k, 0}, \lceil \sigma_{k} \rceil_{n}}$ along $[\alpha_{k, 1}^{n}, \infty)$. Hence the two paths are independent along $[\alpha_{k, 1}^{n}, \infty)$.
\end{proof}
\end{customlemma}
\text{ }
\newline
\begin{customlemma}{5.4} \label{L5.2.4}
 Let $k \in \mathbb{N}$. Let $\Xi_{k}$ be an itinerary of length $k$. Let $\mathcal{B}_{1}, \mathcal{B}_{2}$ be independent standard Brownian motions started at space-time point $(0, \sigma_{k})$. Then as $n \to \infty$:
$$R_{\Xi_{k}}^{n}(t) - R_{\Xi_{k}}^{n}(\lceil \sigma_{k} \rceil_{n}) \overset{d}{\to} \Lambda_{R} (\mathcal{B}_{1}, \mathcal{B}_{2})(t) \text{ in } (\Pi, d) \text{ if } \eta_{k} = 1$$
$$R_{\Xi_{k}}^{n}(t) - R_{\Xi_{k}}^{n}(\lceil \sigma_{k} \rceil_{n}) \overset{d}{\to} \Lambda_{L} (\mathcal{B}_{1}, \mathcal{B}_{2})(t) \text{ in } (\Pi, d) \text{ if } \eta_{k} = -1$$
\begin{proof} $S_{\Xi_{k}}^{n}$ and $Y^{n}_{\beta_{k, 0}^{n}, \lceil \sigma_{k} \rceil_{n}}$ are distributed as random walks with increment $\zeta$ under $(\sqrt{n}, n)$ space-time rescaling. Since $\zeta$ is centered with unit variance, by Donsker's theorem (\cite{DON}) we have:
$$S_{\Xi_{k}}^{n}(t) - S_{\Xi_{k}}^{n}(\lceil \sigma_{k} \rceil_{n}) \overset{d}{\to} \mathcal{B}_{1} \text{ in } (\Pi, d) \text{ as } n \to \infty$$ 
$$Y^{n}_{\beta_{k, 0}^{n}, \lceil \sigma_{k} \rceil_{n}}(t) - Y^{n}_{\beta_{k, 0}^{n}, \lceil \sigma_{k} \rceil_{n}}( \lceil \sigma_{k} \rceil_{n}) \overset{d}{\to} \mathcal{B}_{2} \text{ in } (\Pi, d) \text{ as } n \to \infty$$ 
By Lemma \ref{L5.2.3}, $S_{\Xi_{k}}^{n}$ and $Y^{n}_{\beta_{k, 0}^{n}, \sigma_{k}^{n}}$ are independent along $[\alpha_{k, 1}^{n}, \infty)$, and by definition they are equal to one another along $[\lceil \sigma_{k} \rceil_{n}, \alpha_{k, 1}^{n})$, however this interval shrinks to zero as $n \to \infty$. It follows that:
$$(S_{\Xi_{k}}^{n}(t) - S_{\Xi_{k}}^{n}(\lceil \sigma_{k} \rceil_{n}), Y^{n}_{\beta_{k, 0}^{n}, \lceil \sigma_{k} \rceil_{n}}(t) - Y^{n}_{\beta_{k, 0}^{n}, \lceil \sigma_{k} \rceil_{n}}(\lceil \sigma_{k} \rceil_{n})) \overset{d}{\to} (\mathcal{B}_{1}, \mathcal{B}_{2}) \text{ in } (\Pi, d) \text{ as } n \to \infty$$
By Lemma \ref{L9.0.3}, $\Lambda_{L}$ and $\Lambda_{R}$ are continuous, hence by the continuous mapping theorem we have:
$$\Lambda_{R}(S_{\Xi_{k}}^{n}(t) - S_{\Xi_{k}}^{n}(\lceil \sigma_{k} \rceil_{n}), Y^{n}_{\beta_{k, 0}^{n}, \lceil \sigma_{k} \rceil_{n}}(t) - Y^{n}_{\beta_{k, 0}^{n}, \lceil \sigma_{k} \rceil_{n}}( \lceil \sigma_{k} \rceil_{n})) \overset{d}{\to} \Lambda_{R} (\mathcal{B}_{1}, \mathcal{B}_{2})(t) \text{ in } (\Pi, d)$$
$$\Lambda_{L}(S_{\Xi_{k}}^{n}(t) - S_{\Xi_{k}}^{n}(\lceil \sigma_{k} \rceil_{n}), Y^{n}_{\beta_{k, 0}^{n}, \lceil \sigma_{k} \rceil_{n}}(t) - Y^{n}_{\beta_{k, 0}^{n}, \lceil \sigma_{k} \rceil_{n}}( \lceil \sigma_{k} \rceil_{n})) \overset{d}{\to} \Lambda_{L} (\mathcal{B}_{1}, \mathcal{B}_{2})(t) \text{ in } (\Pi, d)$$
\end{proof}
\end{customlemma}
\subsection{Proof of Theorem 5.6}
Recall that by definition, itineraries $\Xi_{k}$ are composed of rational points and are therefore countable. For each $k \in \mathbb{Z}_{\geq 0}$, fix any enumeration of the pairs $(\Xi_{k, i}, t_{k, i})$ with 
$$\Xi_{k, i} = (x_{i}, s_{i}, (\sigma_{1, i}, ..., \sigma_{k, i}), (\eta_{1, i}, ..., \eta_{k, i}))$$
and such that $t_{k, i} > \sigma_{k, i}$ for $k > 0$ and $t_{k, i}>s_{i}$ for $k=0$. Under these enumerations, let $\mathcal{S}_{k}, \mathcal{S}_{k}^{n} \in \mathbb{R}^{\mathbb{N}}$ be defined as:
$$\mathcal{S}_{k} := \left( \mathcal{G}_{\Xi_{k, 1}}(t_{k,1}), \mathcal{G}_{\Xi_{k, 2}}(t_{k,2}), ... \right)$$
$$\mathcal{S}_{k}^{n} := \left( G_{\Xi_{k, 1}}^{n}(\lceil t_{k,1} \rceil_{n}), G_{\Xi_{k, 2}}^{n}( \lceil t_{k,2} \rceil_{n}), ... \right)$$
Define $\mu_{\infty}$ as the product metric on $\mathbb{R}^{\mathbb{N}}$. That is:
$$\mu_{\infty}\left( (x_{1}, x_{2}, ...), (y_{1}, y_{2}, ... ) \right) := \sum\limits_{i=1}^{\infty} 2^{-i} \frac{|x_{i}-y_{i}|}{1+|x_{i}-y_{i}|}$$
In order to present the proof of Theorem 5.6, we first deal with the case where $k$ is finite:
\begin{customthm}{5.5} \label{T5.3.1} Suppose $\mathcal{V}^{n} \overset{d}{\to} \mathcal{W}$ in $(\mathcal{H}, d_{\mathcal{H}})$. Then for any $k \in \mathbb{Z}_{\geq 0}$ there exists a coupling of $\mathcal{V}^{n}$ and $\mathcal{W}$ such that the following almost surely holds:
$$\left( \mathcal{V}^{n}, \mathcal{S}_{0}^{n}, ..., \mathcal{S}_{k}^{n} \right) \overset{a.s.}{\to} \left( \mathcal{W}, \mathcal{S}_{0}, ..., \mathcal{S}_{k} \right) \text{ in } \left( \mathcal{H} \times \left( \mathbb{R}^{\mathbb{N}} \right)^{k}, d_{\mathcal{H}} \times \mu_{\infty}^{k} \right) \text{ as } n \to \infty$$
\begin{proof} We proceed by induction on $k$.
\newline
\newline
\textbf{Base case:} Note by the Skorokhod representation theorem (\cite{BILL}), there exists a coupling of $\mathcal{V}^{n}$ and $\mathcal{W}$ such that $\mathcal{V}^{n} \overset{a.s.}{\to} \mathcal{W}$. We will show that under such a coupling, $\mathcal{S}_{0}^{n} \overset{a.s.}{\to} \mathcal{S}_{0}$. Fix any pair $(\Xi_{0}, t)$. As argued in the proof of Theorem \ref{T4.0.2}, if $\Xi_{0} = (x,s)$ then:
\begin{equation*}
G_{\Xi_{k}}^{n} = Y_{ \lceil x \rceil_{\sqrt{n}}, \lceil s \rceil_{n}}^{n} \overset{a.s.}{\to} \mathcal{G}_{\Xi_{0}} \text{ in } (\Pi, d)
\end{equation*}
As this argument can be made for any length zero itinerary $\Xi_{0} = (x,s)$ and any $t>s$, this shows coordinate wise convergence of $\mathcal{S}_{0}^{n}$ to $\mathcal{S}_{0}$. Since convergence in $(\mathbb{R}^{\mathbb{N}}, \mu_{\infty})$ is determined by coordinate-wise convergence, it follows that under such a coupling we have:
$$\left( \mathcal{V}^{n}, \mathcal{S}_{0}^{n} \right) \overset{a.s.}{\to} \left( \mathcal{W}, \mathcal{S}_{0} \right) \text{ in } \left( \mathcal{H} \times \mathbb{R}^{\mathbb{N}}, d_{\mathcal{H}} \times \mu_{\infty} \right)$$
\textbf{Inductive step:} Assume that the statement holds for $k \geq 1$. Fix some coupling of $\mathcal{V}^{n}$ and $\mathcal{W}$ such that $\left( \mathcal{V}^{n}, \mathcal{S}_{0}^{n}, ..., \mathcal{S}_{k}^{n} \right) \overset{a.s.}{\to} \left( \mathcal{W}, \mathcal{S}_{0}, ..., \mathcal{S}_{k} \right)$. First we claim that under such a coupling, for any length $k+1$ itinerary, $\Xi_{k+1}$ and any $t > \sigma_{k+1}$, it is the case that:
$$G_{\Xi_{k+1}}^{n}( \lceil t \rceil_{n}) \overset{\mathbb{P}}{\to} \mathcal{G}_{\Xi_{k+1}}(t) \text{ as } n \to \infty$$
Fix a length $k+1$ itinerary $\Xi_{k+1}$ and some $t \in \mathbb{Q}$ with $t>\sigma_{k+1}$. Let $\hat{\Xi}_{k}$ be the length $k$ restriction of $\Xi_{k+1}$. By assumption, $\mathcal{S}_{k}^{n} \overset{a.s.}{\to} \mathcal{S}_{k}$. Since $\sigma_{k+1} \in \mathbb{Q}$, it follows that:
$$G_{\hat{\Xi}_{k}}^{n}(\lceil \sigma_{k+1} \rceil_{n}) \overset{a.s.}{\to} \mathcal{G}_{\hat{\Xi}_{k}}(\sigma_{k+1}) \implies$$

 \begin{equation}\label{E5.3.1}
G_{\Xi_{k+1}}^{n}( \lceil \sigma_{k+1} \rceil_{n}) \overset{a.s.}{\to} \mathcal{G}_{\Xi_{k+1}}(\sigma_{k+1})
\tag{5.5.1}
  \end{equation}
By Lemma \ref{L5.2.4} for $t>\sigma_{k+1}$ we have:
$$R_{\Xi_{k+1}}^{n}(\lceil t \rceil_{n}) - R_{\Xi_{k+1}}^{n}(\lceil \sigma_{k+1} \rceil_{n}) \overset{d}{\to} \Lambda_{R} (\mathcal{B}_{1}, \mathcal{B}_{2})(t)  \text{ if } \eta_{k+1} = 1$$
$$R_{\Xi_{k+1}}^{n}(\lceil t \rceil_{n}) - R_{\Xi_{k+1}}^{n}(\lceil \sigma_{k+1} \rceil_{n}) \overset{d}{\to} \Lambda_{L} (\mathcal{B}_{1}, \mathcal{B}_{2})(t)  \text{ if } \eta_{k+1} = -1$$
Where $\mathcal{B}_{1}, \mathcal{B}_{2}$ are independent standard Brownian motions started at space-time point $(0, \sigma_{k+1})$.
\newline
\newline
By Lemma \ref{L1.3}, the right (left) boundary of all points $(y,t)$ such that $D_{BR}((x, s), (y, t)) \leq 1$ is distributed as a right (left) Skorokhod reflection of two independent Brownian motions started at the point $(x,s)$. Thus:
$$\mathcal{G}_{\Xi_{k+1}}(t)-\mathcal{G}_{\Xi_{k+1}}(\sigma_{k+1}) \overset{d}{=} \Lambda_{R} (\mathcal{B}_{1}, \mathcal{B}_{2})(t) \text{ if } \eta_{k+1} = 1$$
$$\mathcal{G}_{\Xi_{k+1}}(t)-\mathcal{G}_{\Xi_{k+1}}(\sigma_{k+1}) \overset{d}{=} \Lambda_{L} (\mathcal{B}_{1}, \mathcal{B}_{2})(t) \text{ if } \eta_{k+1} = -1
$$
Hence:
 \begin{equation}\label{E5.3.2}
R_{\Xi_{k+1}}^{n}(\lceil t \rceil_{n}) - R_{\Xi_{k+1}}^{n}(\lceil \sigma_{k+1} \rceil_{n}) \overset{d}{\to} \mathcal{G}_{\Xi_{k+1}}(t)-\mathcal{G}_{\Xi_{k+1}}(\sigma_{k+1})
\tag{5.5.2}
\end{equation}
By the Markov property, $R_{\Xi_{k+1}}^{n}(\lceil t \rceil_{n}) - R_{\Xi_{k+1}}^{n}(\lceil \sigma_{k+1} \rceil_{n})$ is independent of $G_{\Xi_{k+1}}^{n}(\lceil \sigma_{k+1} \rceil_{n})$. Similarly, $\mathcal{G}_{\Xi_{k+1}}(t)-\mathcal{G}_{\Xi_{k+1}}(\sigma_{k+1})$ is independent of $\mathcal{G}_{\Xi_{k+1}}(\sigma_{k+1})$. It follows from   \eqref{E5.3.1} and \eqref{E5.3.2} that:
$$( R_{\Xi_{k+1}}^{n}(\lceil t \rceil_{n}) - R_{\Xi_{k+1}}^{n}( \lceil \sigma_{k+1} \rceil_{n}), G_{\Xi_{k+1}}^{n}(\lceil \sigma_{k+1} \rceil_{n}) ) \overset{d}{\to} (\mathcal{G}_{\Xi_{k+1}}(t)-\mathcal{G}_{\Xi_{k+1}}(\sigma_{k+1}), \mathcal{G}_{\Xi_{k+1}}(\sigma_{k+1}) )$$
Hence:
\begin{equation*}
\begin{split}
R^{n}_{\Xi_{k+1}}(\lceil t \rceil_{n}) &= \left[ R_{\Xi_{k+1}}^{n}(\lceil t \rceil_{n}) - R_{\Xi_{k+1}}^{n}(\lceil \sigma_{k+1} \rceil_{n}) \right] + G_{\Xi_{k+1}}^{n}(\lceil \sigma_{k+1} \rceil_{n}) \\
&\overset{d}{\to} \left[ \mathcal{G}_{\Xi_{k+1}}(t)-\mathcal{G}_{\Xi_{k+1}}(\sigma_{k+1}) \right] + \mathcal{G}_{\Xi_{k+1}}(\sigma_{k+1}) \\
&= \mathcal{G}_{\Xi_{k+1}}(t)
\end{split}
\end{equation*}
Without loss of generality we assume $\eta_{k} = 1$. As $\mathcal{V}^{n} \overset{a.s.}{\to} \mathcal{W}$, by Theorem \ref{T4.0.2} and Lemma \ref{L5.2.2} we have (a.s.):
$$\limsup_{n \to \infty} R_{\Xi_{k+1}}^{n}(\lceil t \rceil_{n}) \leq \limsup_{n \to \infty} G_{\Xi_{k+1}}^{n}( \lceil t \rceil_{n}) \leq \mathcal{G}_{\Xi_{k+1}}(t)$$
Combining this with the distributional convergence shown above, by Lemma \ref{L9.0.4} we have:
$$R_{\Xi_{k+1}}^{n}(\lceil t \rceil_{n}) \overset{\mathbb{P}}{\to} \mathcal{G}_{\Xi_{k+1}}(t) \text{ as } n \to \infty$$
Next, define:
$$U^{n} := \left[ G^{n}_{\Xi_{k+1}}(\lceil t \rceil_{n})-\mathcal{G}_{\Xi_{k+1}}(t) \right] \wedge 0$$ 
Then, by Lemma \ref{L5.2.2} we have:
\begin{equation*}
\begin{split}
\mathbb{P} \left( \left| G_{\Xi_{k+1}}^{n}(\lceil t \rceil_{n}) - \mathcal{G}_{\Xi_{k+1}}(t) \right| > \epsilon \right) &= \mathbb{P} \left( U^{n} > \epsilon \right) + \mathbb{P} \left( \mathcal{G}_{\Xi_{k+1}}(t) - G_{\Xi_{k+1}}^{n}(\lceil t \rceil_{n}) > \epsilon \right) \\
&\leq \mathbb{P} \left( U^{n} > \epsilon \right) + \mathbb{P} \left( \mathcal{G}_{\Xi_{k+1}}(t) - R_{\Xi_{k+1}}^{n}(\lceil t \rceil_{n}) > \epsilon \right)
\end{split}
\end{equation*}
Since $\mathcal{V}^{n} \overset{a.s.}{\to} \mathcal{W}$, by Theorem \ref{T4.0.2} we know that $U^{n} \overset{a.s.}{\to} 0$, hence $U^{n} \overset{\mathbb{P}}{\to} 0$. Thus the first term also shrinks to zero, and we have already shown the second term shrinks to zero, thus:
$$G_{\Xi_{k+1}}^{n}(\lceil t \rceil_{n}) \overset{\mathbb{P}}{\to} \mathcal{G}_{\Xi_{k+1}}(t)$$
As the convergence is in probability, the same argument can be made simultaneously for any finite collection of pairs $(\Xi_{k+1}, t)$. Coordinate wise convergence in probability then implies:
$$\mathcal{S}_{k+1}^{n} \overset{\mathbb{P}}{\to} \mathcal{S}_{k+1} \text{ in } (\mathbb{R}^{\mathbb{N}}, \mu_{\infty})$$
Combined with the almost sure convergence given by our inductive assumption, we have:
$$\left( \mathcal{V}^{n}, \mathcal{S}_{0}^{n}, ..., \mathcal{S}_{k+1}^{n} \right) \overset{\mathbb{P}}{\to} \left( \mathcal{W}, \mathcal{S}_{0}, ..., \mathcal{S}_{k+1} \right)$$
Since $\mathcal{H}$ and $\mathbb{R}^{\mathbb{N}}$ are separable, so too is $(\mathcal{H}, \mathbb{R}^{\mathbb{N}}, ..., \mathbb{R}^{\mathbb{N}})$. Therefore we can again apply the Skorokhod representation theorem to find a coupling of $\mathcal{V}^{n}$ and $\mathcal{W}$ such that:
$$\left( \mathcal{V}^{n}, \mathcal{S}_{0}^{n}, ..., \mathcal{S}_{k+1}^{n} \right) \overset{a.s.}{\to} \left( \mathcal{W}, \mathcal{S}_{0}, ..., \mathcal{S}_{k+1} \right)$$
\end{proof}
\end{customthm}
\text{ }
\newline
\noindent We now extend Theorem \ref{T5.3.1} to infinitely many $k$ and present Theorem 5.6. The phrasing of the theorem is slightly different than the one given in the statement of main results but the two are equivalent. Fix any enumeration of \textbf{all} pairs $(\Xi_{k_{i}, i}, t_{i})$ where $\Xi_{k_{i}, i}$ is a length $k_{i}$ itinerary and $t_{i}> \sigma_{k_{i}, i}$ (or $t_{i} > s_{i}$ for $k=0$.) Define:
$$\mathcal{S}_{\infty} := \left( \tanh(\mathcal{G}_{\Xi_{k_{1}, 1}}(t_{1})), \tanh(\mathcal{G}_{\Xi_{k_{2}, 2}}(t_{2})), ... \right)$$
$$\mathcal{S}_{\infty}^{n} := \left( \tanh(G_{\Xi_{k_{1}, 1}}^{n}(\lceil t_{1} \rceil_{n})), \tanh(G_{\Xi_{k_{2}, 2}}^{n}(\lceil t_{2} \rceil_{n})), ... \right)$$

\begin{customthm}{5.6} \label{T5.3.2} Suppose $\mathcal{V}^{n} \overset{d}{\to} \mathcal{W}$ in $(\mathcal{H}, d_{\mathcal{H}})$. Then there exists a coupling of $\mathcal{V}^{n}$ and $\mathcal{W}$ such that the following almost surely holds:
$$\left( \mathcal{V}^{n}, \mathcal{S}_{\infty}^{n} \right) \overset{a.s.}{\to} \left( \mathcal{W}, \mathcal{S}_{\infty} \right) \text{ in } \left( \mathcal{H} \times \mathbb{R}^{\mathbb{N}}, d_{\mathcal{H}} \times \mu_{\infty} \right) \text{ as } n \to \infty$$
\begin{proof} By Theorem \ref{T5.3.1}, for any $i \in \mathbb{N}$ there exists a coupling of $\mathcal{V}^{n}$ and $\mathcal{W}$ such that:
$$\left( \mathcal{V}^{n}, G_{\Xi_{k_{1}, 1}}^{n}(\lceil t_{1} \rceil_{n}), ..., G_{\Xi_{k_{i}, i}}^{n}(\lceil t_{i} \rceil_{n}) \right) \overset{a.s.}{\to} \left( \mathcal{W}, \mathcal{G}_{\Xi_{k_{1}, 1}}(t_{1}), ..., \mathcal{G}_{\Xi_{k_{i}, i}}(t_{i}) \right)$$
It follows that for any $i \in \mathbb{N}$ it is the case that:
$$\left( \mathcal{V}^{n}, G_{\Xi_{k_{1}, 1}}^{n}(\lceil t_{1} \rceil_{n}), ..., G_{\Xi_{k_{i}, i}}^{n}(\lceil t_{i} \rceil_{n}) \right) \overset{d}{\to} \left( \mathcal{W}, \mathcal{G}_{\Xi_{k_{1}, 1}}(t_{1}), ..., \mathcal{G}_{\Xi_{k_{i}, i}}(t_{i}) \right)$$
Since $\mathcal{S}_{\infty}^{n} \in [-1, 1]^{\mathbb{N}}$, for all $n$, the sequence $\{ \mathcal{S}_{\infty}^{n} \}$ is tight. As shown above, the finite-dimensional marginals of $\mathcal{S}^{n}_{\infty}$ converge to those of $\mathcal{S}_{\infty}$, hence by Lemma \ref{L1.6}, we have:
$$\left( \mathcal{V}^{n}, \mathcal{S}_{\infty}^{n} \right) \overset{d}{\to} \left( \mathcal{W}, \mathcal{S}_{\infty} \right)$$
Applying the Skorokhod representation theorem gives the desired result.
\end{proof}
\end{customthm}
\section{Epigraph convergence of $D_{RW}^{n}$ to $D_{BR}$}
We now aim to use Theorem \ref{T5.3.2} in order to prove epigraph convergence. The strategy is to show that any geodesic in the Brownian web can be well approximated by a journey.
\subsection{Approximating geodesics with journeys}
\begin{customthm}{6.1} \label{T6.1.1} The following almost surely holds:
$$\forall (x,s,y,t) \in \mathbb{R}^{4}, \epsilon > 0, $$
$$D_{BR}((x, s), (y,t)) = k < \infty \implies$$
$$\text{There exists an itinerary of length $k$, } \Xi_{k} = (\hat{x}, \hat{s}, \overline{\sigma}, \overline{\eta}) \text{ and } \hat{t} \in \mathbb{Q} \text{ such that:}$$
\begin{itemize}
\begin{multicols}{2}
\item $|\hat{x}-x| < \epsilon$
\item $|\hat{s}-s| < \epsilon$
\item $|\hat{t}-t| < \epsilon$
\item $|\mathcal{G}_{\Xi_{k}}(\hat{t}) - y| < \epsilon$
\end{multicols}
\end{itemize}
\begin{proof} Let $(\gamma_{0}, ..., \gamma_{k})$ be a geodesic between the points $(x,s)$ and $(y,t)$, with $\gamma_{i}$ starting at time $\tau_{i}$. Then:
$$s = \tau_{0} < \tau_{1} < ... < \tau_{k} < t$$
As shown in Theorem 2.1 of \cite{FINR}, $\mathcal{W}$ can be defined as the closure of the set of paths $\{ \gamma^{n} \}$ started at rational points $(x^{n}, s^{n})$. Further, by Lemma \ref{L1.4}, if a sequence of paths $\gamma^{n} \in \mathcal{W}$ converge to $\gamma_{0}$, then the time at which $\gamma^{n}$ and $\gamma_{0}$ coalesce converges to $\tau_{0}$. It follows that for any $\hat{\epsilon}>0$, there exists a path $\hat{\gamma}_{0} \in \mathcal{W}$ starting at space-time point $(\hat{x}, \hat{s}) \in \mathbb{Q}^{2}$ and such that:
\begin{itemize}
\item $|\hat{x}-x| < \epsilon$
\item $|\hat{s}-s| < \epsilon$
\item $\hat{\gamma}_{0}$ coalesces with $\gamma_{0}$ before the time $s+\hat{\epsilon}$
\end{itemize}
If $k=0$ we find such a path $\hat{\gamma}_{0}$ which coalesces with $\gamma_{0}$ before time $\frac{t+\tau_{0}}{2}$, and let $\Xi_{0} = (\hat{x}, \hat{s})$.
\newline
\newline
Suppose now $k>0$. For the sake of notation, we assume $k>1$; the case $k=1$ follows by the same logic. Take some approximating path $\hat{\gamma}$ starting at $(\hat{x}, \hat{s})$ as in the case of $k=0$, but now choose $\hat{\gamma}_{0}$ such that it coalesces with $\gamma_{0}$ before time $\frac{\tau_{0}+\tau_{1}}{2}$. Without loss of generality, we suppose that $\gamma_{1}$ is a right excursion of $\gamma_{0}$. That is to say:
$$\gamma_{1}(\alpha) > \gamma_{0}(\alpha) \text{ for all } \alpha \text{ s.t. } \tau_{1} < \alpha \leq \tau_{2}$$ 
Define the right forward-time distance one boundary of the point $(z, \tau)$ as:
$$\mathfrak{b}_{z, \tau}^{R}(\alpha) := \sup \{ \beta : D_{BR} ((z, \tau), (\beta, \alpha)) \leq 1 \}$$
Define $\lambda := \frac{\tau_{1}+\tau_{2}}{2}$. We claim that $\exists \delta > 0$ such that for all $\alpha$ with $\lambda < \alpha \leq \tau_{2}$ and all $\tau$ with $\tau_{1}-\delta < \tau \leq \tau_{1}$, the following holds:
$$\mathfrak{b}_{\gamma_{0}(\tau), \tau}^{R}(\alpha) = \gamma_{1}(\alpha)$$
Note that since $\gamma_{1}$ is an excursion of $\gamma_{0}$, trivially for $\tau \leq \tau_{1}$ and $\alpha \geq \tau_{1}$ we have:
$$\mathfrak{b}_{\gamma_{0}(\tau), \tau}^{R}(\alpha) \geq \gamma_{1}(\alpha)$$
By way of contradiction, suppose for all $n \in \mathbb{N}$ there exists $\tau^{n}$ with $\tau_{1} - \frac{1}{n} < \tau^{n} \leq \tau_{1}$ and $\alpha^{n}$ with $\lambda < \alpha^{n} \leq \tau_{2}$ such that:
$$\mathfrak{b}_{\gamma_{0}(\tau^{n}), \tau^{n}}^{R}(\alpha^{n}) > \gamma_{1}(\alpha^{n})$$
First we claim this would imply there exists a sequence of right excursions, $\gamma^{n}$ started at time $t^{n} \in [\tau^{n}, \tau_{1}]$ such that $\gamma^{n}(\alpha^{n}) > \gamma_{0}(\alpha^{n})$. If not, there would instead need to exist right excursions started at times $t^{n}$ with $\tau_{1} < t^{n} \leq \tau_{2}$ which surpass $\gamma_{1}$. Since $\gamma_{0}$ and $\gamma_{1}$ fall along a geodesic between $(x,s)$ and $(y,t)$, we know for $w$ such that $\tau_{1} < w \leq \tau_{2}$ it must be the case that $\gamma_{0}(w) < \gamma_{1}(w)$. Due to the coalesence property of $\mathcal{W}$, this implies any right excursions of $\gamma_{0}$ started at such times could not surpass $\gamma_{1}$, else they would coalesce with it. Therefore there must  exist a sequence of right excursions $\gamma^{n}$ started at times $t^{n}$ with $\tau_{1} -\frac{1}{n} < t^{n} \leq \tau_{1}$. Note then that:
\begin{itemize}
\item $\tau_{1} - \frac{1}{n} \leq t^{n} \leq \tau_{1} \implies \lim_{n \to \infty} t^{n} = \tau_{1}$
\item $\gamma^{n}(\alpha^{n}) > \gamma_{1}(\alpha^{n})$
\end{itemize}
By the compactness of $\mathcal{W}$, there would exist a subsequential limit $\gamma^{\infty} := \lim_{j \to \infty} \gamma^{n_{j}}$. By the fact that $\gamma_{0} \in \Pi^{\star}$ we would have:
$$\lim_{j \to \infty} (\gamma_{0}(t^{n_{j}}), t^{n_{j}}) = (\gamma_{0}(\tau_{1}), \tau_{1})$$
Thus $\gamma^{\infty}$ would be a path starting at the point $(\gamma_{0}(\tau_{1}), \tau_{1})$. We claim $\gamma^{\infty}$ would be distinct from $\gamma_{1}$ and from the continuation of $\gamma_{0}$. Note that:
$$\gamma^{n}(\alpha^{n}) > \gamma_{1}(\alpha^{n}) > \gamma_{0}(\alpha^{n})$$
Thus it cannot be the case that $\gamma^{n}$ coalesces with $\gamma_{0}$ or $\gamma_{1}$ before the time $\lambda$. Then by Lemma \ref{L1.4}, it cannot be the case that $\gamma^{n}$ converges to the excursion $\gamma_{1}$, nor to the continuation of $\gamma_{0}$ started at the time $\tau_{1}$. This would imply that there exist three distinct outgoing paths from $(\gamma_{0}(\tau_{1}), \tau_{1})$. There is also at least one incoming path - $\gamma_{0}$. This would contradict Lemma \ref{L1.5}, hence this cannot be true.
\newline
\newline
Therefore no such sequence $\gamma^{n}$ can exist, and thus there must exist a $\delta > 0$ such that:
$$\tau_{1} - \delta < \tau \leq \tau_{1} \implies$$
$$\mathfrak{b}_{\gamma_{0}(\tau), \tau}^{R}(\alpha) = \gamma_{1}(\alpha) \text{ for all } \alpha \text{ s.t. } \lambda < \alpha \leq \tau_{2}$$
Find any such $\delta$. Take any $\sigma_{1} \in \mathbb{Q}$ with $\left( \tau_{1} - \delta \vee \frac{\tau_{0}+\tau_{1}}{2} \right) < \sigma_{1} \leq \tau_{1}$. Define $\Xi_{1} := (\hat{x}, \hat{s}, \sigma_{1}, 1)$. Since $\hat{\gamma}_{0}$ coalesces with $\gamma_{0}$ before time $\frac{\tau_{0}+\tau_{1}}{2}$, for $\alpha$ with $\lambda < \alpha \leq \tau_{2}$ we have:
\begin{equation*}
\begin{split}
\mathcal{G}_{\Xi_{1}}(\alpha) &= \mathfrak{b}_{\hat{\gamma_{0}}(\sigma_{1}), \sigma_{1}}^{R}(\alpha) \\
&=\mathfrak{b}_{\gamma_{0}(\sigma_{1}), \sigma_{1}}^{R}(\alpha) \\
&= \gamma_{1}(\alpha)
\end{split}
\end{equation*}
By iterating this process, for any $k \in \mathbb{N}$ we can find an itinerary $\Xi_{k} := (\hat{x}, \hat{s}, \overline{\sigma}, \overline{\eta})$ such that for all $\alpha$ with $\frac{\tau_{k}+t}{2} < \alpha \leq t$ it is the case that:
$$\mathcal{G}_{\Xi_{k}}(\alpha) = \gamma_{k}(\alpha)$$
By the continuity of $\gamma_{k}$, we can then find $\hat{t} \in \mathbb{Q}$ with $\frac{\tau_{k}+t}{2} < \hat{t} \leq t$ such that:
\begin{itemize}
\item $|\hat{t} - t| < \epsilon$
\item $|\gamma_{k}(\hat{t})-\gamma_{k}(t)| < \epsilon$
\end{itemize}
This gives the desired pair, $(\Xi_{k}, \hat{t})$.
\end{proof}
\end{customthm}
\begin{figure}[h] \label{Fig5}
    \centering
    \includegraphics[width=13cm]{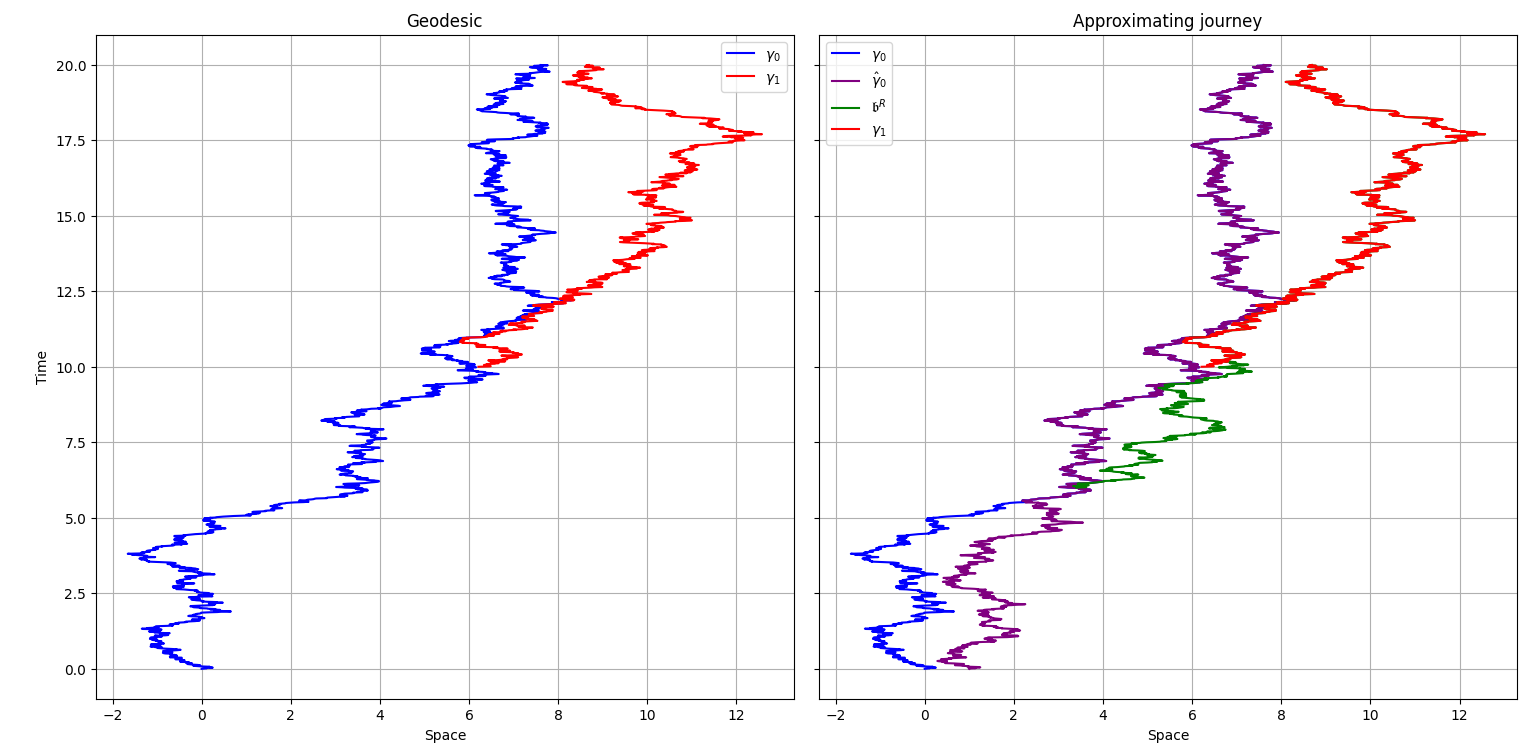}
    \caption{$\hat{\gamma}_{0}$, a path chosen with starting points $(\hat{x}, \hat{s}) \in \mathbb{Q}^{2}$, can be chosen as so that $\hat{\gamma}_{0}$ and $\gamma_{0}$ coalesce prior to the mid-point between $\tau_{0}$ and $\tau_{1}$. $\mathfrak{b}^{R}$ is the right boundary of all points of distance one from a point $(\gamma_{0}(\tau), \tau)$. As $\tau> \frac{\tau_{0}+\tau_{1}}{2}$, it is also the right boundary of all points of distance one from the point $(\hat{\gamma}_{0}(\tau), \tau)$. If $\tau$ is chosen close enough to $\tau_{1}$, we can be guaranteed that $\mathfrak{b}^{R}$ and $\gamma_{1}$ will coalesce before the mid-point between $\tau_{1}$ and $\tau_{2}$. Note the term coalesce is used loosely here as $\mathfrak{b}^{R}$ is not an element of $\mathcal{W}$; after $\gamma_{0}$ and $\gamma_{1}$ have coalesced, $\mathfrak{b}^{R}$ may no longer agree with $\gamma_{1}$. Fortunately due to the fact that $\gamma_{0}, \gamma_{1}$ comprise a geodesic, $\gamma_{0}$ and $\gamma_{1}$ cannot coalesce before the time $\tau_{2}$. The fact that $\mathfrak{b}^{R}$ and $\gamma_{1}$ agree allows us iterate the process. Above is a depiction of the phenomenon. The original geodesic follows the blue path (which coalesces with the purple path) and then the red path (which coalesces with the green path). The approximating journey follows the purple path and then the green path.}
    \label{fig:geodapprox}
\end{figure}
\subsection{Proof of Theorem 6.3}
We approach Theorem 6.3 by utilizing Lemma \ref{4.1VV}. We use Theorem \ref{T6.1.1} to show that the second condition in Lemma \ref{4.1VV} is satisfied:
\begin{customlemma}{6.2} \label{L6.2.1} Under the coupling given in Theorem \ref{T5.3.2}, the following (almost surely) holds:
$$\forall (x, s, y, t) \in \mathbb{R}^{4}, \text{ there exists a sequence } (x^{n}, s^{n}, y^{n}, t^{n}) \to (x,s,y,t) \text{ such that: }$$
$$\limsup_{n \to \infty} D_{RW}^{n} ((x^{n}, s^{n}), (y^{n}, t^{n})) \leq D_{BR}((x,s),(y,t))$$
\begin{proof} If $D_{BR}((x,s),(y,t)) = \infty$ the result is trivial, so we assume $D_{BR}((x,s),(y,t)) = k \in \mathbb{Z}_{\geq 0}$. By Theorem \ref{T6.1.1}, for $j \in \mathbb{N}$, we can find a sequence of pairs $(\Xi_{k}^{j}, \hat{t}_{j})$ with $\Xi_{k}^{j} = (\hat{x}_{j}, \hat{s}_{j}, \overline{\sigma}_{j}, \overline{\eta}_{j})$ such that:
\begin{itemize}
\begin{multicols}{2}
\item $|\hat{x}_{j}-x| < \frac{1}{j}$
\item $|\hat{s}_{j}-s| < \frac{1}{j}$
\item $|\mathcal{G}_{\Xi_{k}^{j}}(\hat{t}_{j}) - y| < \frac{1}{j}$
\item $|\hat{t}_{j} - t| < \frac{1}{j}$
\end{multicols}
\end{itemize}
Define $\{ n_{j} \} \subseteq \mathbb{N}$ with $n_{0}=1$ and, for $j \geq 1$:
$$n_{j} := \min \{ n : n> n_{j-1}, \left| G_{\Xi_{k}^{j}}^{m}( \lceil \hat{t}_{j} \rceil_{m}) - \mathcal{G}_{\Xi_{k}^{j}}(\hat{t}_{j}) \right| < \frac{1}{j} \text{ for all } m \geq n \}$$
For $n < n_{1}$ let $(x^{n}, s^{n}, y^{n}, t^{n})=(1,1,1,1)$. For $j \geq 1$ and $n_{j} \leq n < n_{j+1}$, define:
$$(x^{n}, s^{n}, y^{n}, t^{n}) := ( \lceil \hat{x}_{j}\rceil_{\sqrt{n}}, \lceil \hat{t}_{j} \rceil_{n}, G_{\Xi_{k}^{j}}^{n}( \lceil \hat{t}_{j} \rceil_{n}), \lceil \hat{t}_{j} \rceil_{n})$$
By the conclusions of Theorem \ref{T5.3.2}, it is the case that $n_{j} < \infty$ for all $j$, and it is then clear that $(x^{n}, s^{n}, y^{n}, t^{n}) \to (x,s,y,t)$. Then we have:
$$\limsup_{n \to \infty} D_{RW}^{n} ((x^{n}, s^{n}), (y^{n}, t^{n})) \leq k = D_{BR}((x,s),(y,t))$$
\end{proof}
\end{customlemma}
\text{ }
\newline
After our preparation, Theorem 6.3 is immediate:
\begin{customthm}{6.3} \label{T6.2.2} Under the coupling given in Theorem \ref{T5.3.2}, the following (almost surely) holds:
$$\mathfrak{e}\left( D_{RW}^{n} \right) \to \mathfrak{e} \left( D_{BR} \right) \text{ in } (\mathcal{E}_{\star}, d_{\star}) \text{ as } n \to \infty$$
\begin{proof} Since $D_{RW}^{n}(u,v) = \infty$ if $u$ or $v$ are not in the rescaled lattice, the functions $D_{RW}^{n}$ are trivially lower semicontinuous. By Lemma \ref{L1.7}, Theorem \ref{T4.0.1}, Lemma \ref{L6.2.1} and Lemma \ref{4.1VV}, the result then follows.
\end{proof}
\end{customthm}
\section{Convergence of journeys in $(\Pi, d)$}
We would like to extend the results of Theorem \ref{T5.3.2} to convergence in $(\Pi, d)$. The issue is that Theorem \ref{T5.3.2} gives pointwise convergence, rather than the uniform convergence on compact sets that characterizes convergence in $(\Pi, d)$. In order to address the issue, we aim to show that $G_{\Xi_{k}}^{n}$ are unlikely to vary too wildly for large enough $n$. As in proving Theorem \ref{T5.3.2}, we approach the problem by first analyzing the behavior of the approximation curves $R_{\Xi_{k}}^{n}$.
\subsection{Definitions relating to continuity and variation}
For $\delta, M > 0$ we define a modulus of continuity function, $\mathcal{M}_{\delta, M} : \Pi^{\star} \to \overline{\mathbb{R}}$:
$$\mathcal{M}_{\delta, M}(f) := \sup_{\substack{|w|, |s| \leq M \\ |w-s| < \delta}} |\hat{f}(w)-\hat{f}(s)|$$
For $w, s \in \mathbb{R}$ with $w < s$, we define the positive and negative total variation functions, $P_{w,s} : \Pi^{\star} \to \overline{\mathbb{R}}$, and $N_{w,s} : \Pi^{\star} \to \overline{\mathbb{R}}$:
$$P_{w,s}(f) := \sup \left\{ \sum\limits_{j=1}^{j=m} \delta_{j}(\hat{f}(y_{j})-\hat{f}(y_{j-1})) : w=y_{0} < y_{1} < ... < y_{m} = s, \delta_{j} \in \{ 0, 1 \} \text{ for all } j \right\}$$
$$N_{w,s}(f) := \sup \left\{ \sum\limits_{j=1}^{j=m} \delta_{j}(\hat{f}(y_{j-1})-\hat{f}(y_{j})) : w=y_{0} < y_{1} < ... < y_{m} = s, \delta_{j} \in \{ 0, 1 \} \text{ for all } j \right\}$$
We will define extensions of our approximation curves, $\mathcal{R}_{\Xi_{k}}^{n}$. For an itinerary $\Xi_{k} = (x, s, \overline{\sigma}, \overline{\eta})$, let $(\mathcal{R}_{\Xi_{k}}^{n}, \lceil s \rceil_{n}) \in \Pi^{\star}$ be defined as:
$$\mathcal{R}_{\Xi_{k}}^{n}(t) := \begin{cases} G_{\Xi_{k}}^{n}(t) & \text{ for } t \in [ \lceil s \rceil_{n}, \lceil \sigma_{k} \rceil_{n}] \\ R_{\Xi_{k}}^{n}(t) & \text{ for } t > \lceil \sigma_{k} \rceil_{n} \end{cases}$$
We define the push processes, used to measure the effect of the map $\Lambda_{R}$ or $\Lambda_{L}$ had on the path $S_{\Xi_{k}}^{n}$. $I_{\Xi_{k}}^{n} : \lceil \sigma_{k} \rceil_{n} \to \mathbb{R}$:
$$I_{\Xi_{k}}^{n}(t) := \left| R_{\Xi_{k}}^{n}(t) - S_{\Xi_{k}}^{n}(t) \right|$$
Lastly we define the error function, measuring the error of our approximation functions. $E_{\Xi_{k}}^{n} : [\lceil s \rceil_{n}, \infty) \to \mathbb{R}$:
$$E_{\Xi_{k}}^{n}(t) :=  \left| G_{\Xi_{k}}^{n}(t) - \mathcal{R}_{\Xi_{k}}^{n}(t) \right|$$

\begin{figure}[h] \label{Fig10}
    \centering
    \includegraphics[width=15cm]{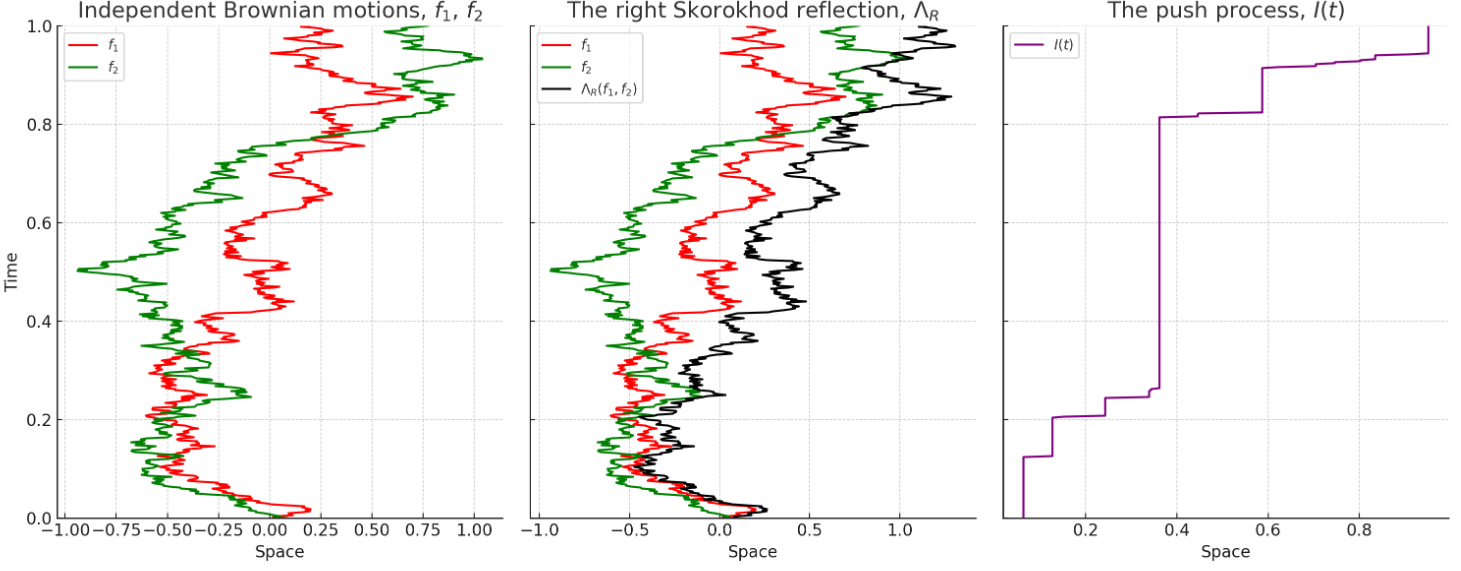}
    \caption{The right Skorokhod reflection $\Lambda_{R}$ and the associated push process $I := \Lambda_{R}(f_{1}, f_{2}) - f_{1}$. Note $I$ is non-decreasing.}
    \label{fig:pushp}
\end{figure}

\subsection{Basic results relating to continuity and variation}

\begin{customlemma}{7.1} \label{L7.2.1}
For any fixed $M, \delta, \epsilon >0$ it is the case that:
$$\lim\limits_{\delta \to 0^{+}} \limsup\limits_{n \to \infty} \mathbb{P}\left( \mathcal{M}_{\delta, M}(I_{\Xi_{k}}^{n}) > \epsilon \right) = 0$$
\begin{proof} Without loss of generality, we suppose $\eta_{k}=1$. By Lemma \ref{L9.0.5}, we have:
 \begin{equation}\label{E7.2.1.1}
\begin{split}
\mathbb{P} \left( \mathcal{M}_{\delta, M}(I_{\Xi_{k}}^{n}) > 4\epsilon \right) &\leq \mathbb{P}\left( \mathcal{M}_{\delta, M}(R_{\Xi_{k}}^{n}) > 3 \epsilon \right) + \mathbb{P}\left( \mathcal{M}_{\delta, M}(S_{\Xi_{k}}^{n}) >  \epsilon \right) \\
&\leq \mathbb{P} \left( \mathcal{M}_{\delta, M} (Y^{n}_{\beta_{k, 0}^{n}, \sigma_{k}^{n}}) > \epsilon \right) + 2 \mathbb{P} \left( \mathcal{M}_{\delta, M}(S_{\Xi_{k}}^{n}) > \epsilon \right) \\
\end{split}
\tag{7.1.1}
\end{equation}
Letting $\mathcal{B}$ represent standard Brownian motion started at $(0, \sigma_{k})$, by Donsker's theorem (\cite{DON}), Lemma \ref{L9.0.6}, and the continuous mapping theorem we have:
$$\mathcal{M}_{\delta, M}(S_{\Xi_{k}}^{n}) \overset{d}{\to} \mathcal{M}_{\delta, M}(\mathcal{B})$$
$$\mathcal{M}_{\delta, M}(Y^{n}_{\beta^{n}_{k, 0}, \lceil \sigma_{k} \rceil_{n}}) \overset{d}{\to} \mathcal{M}_{\delta, M}(\mathcal{B})$$
By L\'evy's modulus of continuity (\cite{LEV}), it follows that:
$$\lim_{\delta \to 0^{+}} \limsup_{n \to \infty} \mathbb{P} \left( \mathcal{M}_{\delta, M}(Y^{n}_{\beta_{k, 0}^{n}, \lceil \sigma_{k} \rceil_{n}}) > \epsilon \right) = 0$$
$$\lim_{\delta \to 0^{+}} \limsup_{n \to \infty} \mathbb{P} \left( \mathcal{M}_{\delta, M}(S_{\Xi_{k}}^{n}) > \epsilon \right) = 0$$
These facts and \eqref{E7.2.1.1} give the desired result.
\newline
\newline
\end{proof}
\end{customlemma}

\begin{customlemma}{7.2} \label{L7.2.3} Let $k \in \mathbb{N}$. Let $\Xi_{k} = (x,s,\overline{\sigma}, \overline{\eta})$. Let $M, \delta >0$ and suppose that $\lceil \sigma_{k} \rceil_{n} \leq w \leq t \leq M$ and that $|t-w| < \delta$. Then:
$$\sup\limits_{a \in [w, t]} E_{\Xi_{k}}^{n}(a) \leq E_{\Xi_{k}}^{n}(t) + \mathcal{M}_{\delta, M}(I_{\Xi_{k}}^{n})$$
\begin{proof} Without loss of generality, suppose $\eta_{k} = 1$. By Lemma \ref{L5.2.1} we have:
\begin{equation} \label{7.2.3.1}
N_{w,t}(E_{\Xi_{k}}^{n}) \leq P_{w,t}(I_{\Xi_{k}}^{n})
\tag{7.2.1}
\end{equation} 
Further, by definition of negative and positive variation, we have:
\begin{equation} \label{7.2.3.2}
\sup\limits_{a \in [w, t]} E_{\Xi_{k}}^{n}(a) \leq N_{w, t}(E_{\Xi_{k}}^{n}) + E_{\Xi_{k}}^{n}(t)
\tag{7.2.2}
\end{equation}
Since $I_{\Xi_{k}}^{n}$ is non-decreasing, we know $P_{w,t}(I_{\Xi_{k}}^{n}) = I_{\Xi_{k}}^{n}(t)-I_{\Xi_{k}}^{n}(w)$. Combining this with \eqref{7.2.3.1} and \eqref{7.2.3.2} we have:
\begin{equation*}
\begin{split}
\sup\nolimits_{a \in [w, t]} E_{\Xi_{k}}^{n}(a) & \leq N_{w,t}(E_{\Xi_{k}}^{n}) + E_{\Xi_{k}}^{n}(t) \\
&\leq P_{w,t}(I_{\Xi_{k}}^{n}) + E_{\Xi_{k}}^{n}(t) \\
& = I_{\Xi_{k}}^{n}(t) - I_{\Xi_{k}}^{n}(w) + E_{\Xi_{k}}^{n}(t) \\
& \leq \mathcal{M}_{\delta, M}(I_{\Xi_{k}}^{n}) + E_{\Xi_{k}}^{n}(t)
\end{split}
\end{equation*}
\end{proof}
\end{customlemma}
\text{ }
\newline

\begin{customlemma}{7.3} \label{L7.2.4} Let $k \in \mathbb{N}$ and let $\Xi_{k} = (x,s,\overline{\sigma}, \overline{\eta})$ be a fixed length $k$ itinerary. Suppose the following conditions hold:
\begin{itemize}
\item $\mathcal{R}_{\Xi_{k}}^{n} \overset{\mathbb{P}}{\to} \mathcal{G}_{\Xi_{k}} \text{ in } (\Pi, d)$ as $n \to \infty$.
\item $\mathcal{V}^{n} \overset{a.s.}{\to} \mathcal{W} \text{ in } (\mathcal{H}, d_{\mathcal{H}})$ as $n \to \infty$.
\item $G_{\Xi_{k}}^{n} (\lceil \sigma_{k} \rceil_{n}) \overset{a.s.}{\to} \mathcal{G}_{\Xi_{k}}(\sigma_{k})$
\end{itemize}
Then for any $t > s , \epsilon >0$ it is the case that $\lim\limits_{n \to \infty} \mathbb{P} \left( E_{\Xi_{k}}^{n}(t) > \epsilon \right) = 0$.
\begin{proof} Without loss of generality, we suppose $\eta_{k} = 1$. If $t \leq \sigma_{k}$ then the statement is trivial as $E_{\Xi_{k}}^{n}(t) = 0$. So we assume $t>\sigma_{k}$, then:
\begin{equation*} \label{E7.2.4.1}
\begin{split}
\mathbb{P} \left(  E_{\Xi_{k}}^{n}(t)  > \epsilon \right) &=\mathbb{P} \left( G_{\Xi_{k}}^{n}(t) - \mathcal{R}_{\Xi_{k}}^{n}(t)  > \epsilon \right) \\
&=\mathbb{P} \left( G_{\Xi_{k}}^{n}(t) - \mathcal{R}_{\Xi_{k}}^{n}(t)  > \epsilon, G_{\Xi_{k}}^{n}(t) > \mathcal{G}_{\Xi_{k}}(t) + \epsilon \right) 
\\
&+ \mathbb{P} \left( G_{\Xi_{k}}^{n}(t) - \mathcal{R}_{\Xi_{k}}^{n}(t)  > \epsilon, G_{\Xi_{k}}^{n}(t) \leq \mathcal{G}_{\Xi_{k}}(t) + \epsilon \right) \\
&\leq \mathbb{P} \left( G_{\Xi_{k}}^{n}(t) > \mathcal{G}_{\Xi_{k}}(t) + \epsilon \right) + \mathbb{P} \left( \mathcal{G}_{\Xi_{k}}(t) - \mathcal{R}^{n}_{\Xi_{k}}(t) > 2\epsilon \right)
\end{split} \tag{7.3.1}
\end{equation*}
By Theorem \ref{T4.0.2}, $\limsup_{n \to \infty} G_{\Xi_{k}}^{n}(t) \overset{a.s.}{\leq} \mathcal{G}_{\Xi_{k}}(t)$, hence the first term on the last line of \eqref{E7.2.4.1} shrinks to zero as $n \to \infty$, and the second term shrinks to zero by our first assumption.
\end{proof}
\end{customlemma}
\text{ }
\newline

\subsection{Proof of Theorem 7.5}
As was the case in proving Theorem \ref{T5.3.2}, we define sequence spaces and first prove the result for finite $k$. In this case we define sequences $\hat{\mathcal{S}} \in \Pi^{\mathbb{N}}$. For each $k \in \mathbb{Z}_{\geq 0}$, fix any enumeration of the itineraries $(\Xi_{k})$. Under these enumerations, let $\hat{\mathcal{S}}_{k}, \hat{\mathcal{S}}_{k}^{n} \in \Pi^{\mathbb{N}}$ be defined as:
$$\hat{\mathcal{S}}_{k} := \left( \mathcal{G}_{\Xi_{k, 1}}, \mathcal{G}_{\Xi_{k, 2}}, ... \right)$$
$$\hat{\mathcal{S}}_{k}^{n} := \left( G_{\Xi_{k, 1}}^{n}, G_{\Xi_{k, 2}}^{n}, ... \right)$$
Define $\hat{\mu}_{\infty}$ as the product metric on $\Pi^{\mathbb{N}}$. That is:
$$\hat{\mu}_{\infty}\left( (x_{1}, x_{2}, ...), (y_{1}, y_{2}, ... ) \right) := \sum\limits_{i=1}^{\infty} 2^{-i} \frac{d(x_{i}, y_{i})}{1+d(x_{i}, y_{i})}$$
\begin{customthm}{7.4} \label{T7.3.1} Suppose $\mathcal{V}^{n} \overset{d}{\to} \mathcal{W}$ in $(\mathcal{H}, d_{\mathcal{H}})$. Then for any $k \in \mathbb{Z}_{\geq 0}$ there exists a coupling of $\mathcal{V}^{n}$ and $\mathcal{W}$ such that the following holds:
$$\left( \mathcal{V}^{n}, \hat{\mathcal{S}}_{0}^{n}, ..., \hat{\mathcal{S}}_{k}^{n} \right) \overset{a.s.}{\to} \left( \mathcal{W}, \hat{\mathcal{S}}_{0}, ..., \hat{\mathcal{S}}_{k} \right) \text{ in } \left( \mathcal{H} \times \left( \Pi^{\mathbb{N}} \right)^{k}, d_{\mathcal{H}} \times \hat{\mu}_{\infty}^{k} \right) \text{ as } n \to \infty$$
\begin{proof} We proceed by induction on $k$.
\newline
\newline
\textbf{Base case:} This follows from the same argument used in the base case of Theorem \ref{T5.3.1}.
\newline
\newline
\textbf{Inductive step:} Let $k \geq 1$. Let $\Xi_{k}$ be a length k itinerary and take some coupling of $\mathcal{V}^{n}$ and $\mathcal{W}$ such that $\mathcal{V}^{n} \overset{a.s.}{\to} \mathcal{W}$ and $\hat{\mathcal{S}}_{k-1}^{n} \overset{a.s.}{\to} \hat{\mathcal{S}}_{k-1}$.
\newline
\newline
We claim that under such a coupling, $\mathcal{R}_{\Xi_{k}}^{n} \overset{\mathbb{P}}{\to} \mathcal{G}_{\Xi_{k}}$. Let $\hat{\Xi}_{k-1}$ be the length $k-1$ restriction of $\Xi_{k}$. Without loss of generality, assume $\eta_{k}=1$. Let $\mathcal{G}_{\Xi_{k}}^{\star}$ be the restriction of $\mathcal{G}_{\Xi_{k}}$ to $[\sigma_{k}, \infty)$. Let $\mathcal{B}_{1}, \mathcal{B}_{2}$ be standard Brownian motions started at $(\mathcal{G}_{\Xi_{k}}(\sigma_{k}), \sigma_{k})$, then by our chosen coupling and Lemma \ref{L5.2.4}, we have:
$$G_{\hat{\Xi}_{k-1}}^{n} \overset{a.s.}{\to} \mathcal{G}_{\hat{\Xi}_{k-1}} \text{ in } (\Pi, d)$$
$$R_{\Xi_{k}}^{n} \overset{d}{\to} \Lambda_{R}\left( \mathcal{B}_{1}, \mathcal{B}_{2} \right) \text{ in } (\Pi, d)$$
By Lemma \ref{L1.3}, $\Lambda_{R}(\mathcal{B}_{1}, \mathcal{B}_{2}) \overset{d}{=} \mathcal{G}_{\Xi_{k}}^{\star}$, hence $\mathcal{R}_{\Xi_{k}}^{n} \overset{d}{\to} \mathcal{G}_{\Xi_{k}}$ in $(\Pi, d)$.
\newline
\newline
By the convergence of $G_{\hat{\Xi}_{k-1}}^{n}$ to $\mathcal{G}_{\hat{\Xi}_{k-1}}$ we have:
$$\limsup_{n \to \infty} \mathcal{R}_{\Xi_{k}}^{n}(t) \overset{a.s.}{=} \mathcal{G}_{\Xi_{k}}(t) \text{ for } t \leq \sigma_{k}$$
While by Theorem \ref{T4.0.2} and Lemma \ref{L5.2.2} we have:
$$\limsup_{n \to \infty} \mathcal{R}_{\Xi_{k}}^{n}(t) \leq \limsup_{n \to \infty} G_{\Xi_{k}}^{n}(t) \overset{a.s.}{\leq} \mathcal{G}_{\Xi_{k}}(t) \text{ for } t > \sigma_{k}$$
Hence by Lemma \ref{L9.0.7}, it follows that $\mathcal{R}_{\Xi_{k}}^{n} \overset{\mathbb{P}}{\to} \mathcal{G}_{\Xi_{k}}$ in $(\Pi, d)$. We claim that $G_{\Xi_{k}}^{n} \overset{\mathbb{P}}{\to} \mathcal{G}_{\Xi_{k}}$. Note that:
\begin{equation*}
\begin{split}
\mathbb{P} \left( d(G_{\Xi_{k}}^{n}, \mathcal{G}_{\Xi_{k}}) > 4 \epsilon \right) &\leq \mathbb{P} \left( d(G_{\Xi_{k}}^{n}, \mathcal{R}_{\Xi_{k}}^{n}) > 2 \epsilon \right) \\
&+ \mathbb{P} \left( d(\mathcal{R}_{\Xi_{k}}^{n}, \mathcal{G}_{\Xi_{k}}) > 2 \epsilon \right)
\end{split}
\end{equation*}
It then suffices to show that $\lim_{n \to \infty} \mathbb{P} \left( d(G_{\Xi_{k}}^{n}, \mathcal{R}_{\Xi_{k}}^{n}) > 2 \epsilon \right)  = 0$. Fix any $\hat{\epsilon} > 0$. Take $M > 0$ such that $(M+1)^{-1} < \epsilon$. By Lemma \ref{L7.2.1} we can find $\delta > 0$ such that:
$$\limsup\limits_{n \to \infty} \mathbb{P} \left( \mathcal{M}_{\delta, M}( I_{\Xi_{k}}^{n} ) > \epsilon \right) < \hat{\epsilon}$$
Take some such $\delta$. Let $-M = t_{0} < t_{1} < ... < t_{i} = M$ be a partition of $[-M, M]$ such that $t_{j}-t_{j-1} < \delta$ for all $1 \leq j \leq i$. By Lemma \ref{L7.2.3} we have:
\begin{equation*}
\begin{split}
\mathbb{P} \left( d(G_{\Xi_{k}}^{n}, \mathcal{R}_{\Xi_{k}}^{n}) > 2 \epsilon \right) &\leq \mathbb{P} \left( \sup\limits_{-M \leq t \leq M} E_{\Xi_{k}}^{n}(t) > 2 \epsilon \right) \\
& \leq \left( \mathbb{P} \left( \mathcal{M}_{\delta, M} (I_{\Xi_{k}}^{n} ) > \epsilon \right) \right) + \left( \sum\limits_{j=1}^{j=i} \mathbb{P} \left( E_{\Xi_{k}}^{n}(x_{j}) > \epsilon \right) \right)
\end{split}
\end{equation*}
By our choice of $\delta$ the first term shrinks to $\hat{\epsilon}$. By Lemma \ref{L7.2.4}, the sum in the second term shrinks to zero. Hence showing for any $\hat{\epsilon} >0$, it is the case that:
$$\lim_{n \to \infty} \mathbb{P} \left( d(G_{\Xi_{k}}^{n}, \mathcal{R}_{\Xi_{k}}^{n}) > 2 \epsilon \right) < \hat{\epsilon}$$
Hence allowing us to conclude that:
$$G_{\Xi_{k}}^{n} \overset{\mathbb{P}}{\to} \mathcal{G}_{\Xi_{k}}$$ 
Since the same argument can be simultaneously made for any finite number of length $k$ itineraries, it follows that:
$$\hat{\mathcal{S}}_{k}^{n} \overset{\mathbb{P}}{\to} \hat{\mathcal{S}}_{k}$$
Therefore there exists a coupling of $\mathcal{V}^{n}$ and $\mathcal{W}$ such that:
$$\left( \mathcal{V}^{n}, \hat{\mathcal{S}}_{0}^{n}, ..., \hat{\mathcal{S}}_{k}^{n} \right) \overset{\mathbb{P}}{\to} \left( \mathcal{W}, \hat{\mathcal{S}}_{0}, ..., \hat{\mathcal{S}}_{k} \right)$$
As shown in \cite{FINR}, $(\Pi, d)$ is separable, hence as is $(\Pi^{\mathbb{N}}, \hat{\mu}_{\infty})$, thus we can apply the Skorokhod representation theorem and argue there exists a coupling of $\mathcal{V}^{n}$ and $\mathcal{W}$ such that:
$$\left( \mathcal{V}^{n}, \hat{\mathcal{S}}_{0}^{n}, ..., \hat{\mathcal{S}}_{k}^{n} \right) \overset{a.s.}{\to} \left( \mathcal{W}, \hat{\mathcal{S}}_{0}, ..., \hat{\mathcal{S}}_{k} \right)$$
Completing the induction.
\end{proof}
\end{customthm}
\text{ }
\newline
\noindent As was done in the proof of Theorem \ref{T5.3.2}, we now extend Theorem \ref{T7.3.1} to infinitely many $k$. Fix any enumeration of \textbf{all} itineraries $\Xi_{k_{i}, i}$. Define:
$$\hat{\mathcal{S}}_{\infty} := \left( \mathcal{G}_{\Xi_{k_{1}, 1}}, \mathcal{G}_{\Xi_{k_{2}, 2}}, ... \right)$$
$$\hat{\mathcal{S}}_{\infty}^{n} := \left( G_{\Xi_{k_{1}, 1}}^{n}, G_{\Xi_{k_{2}, 2}}^{n}, ... \right)$$

\begin{customthm}{7.5} \label{T7.3.2} Suppose $\mathcal{V}^{n} \overset{d}{\to} \mathcal{W}$ in $(\mathcal{H}, d_{\mathcal{H}})$. Then there exists a coupling of $\mathcal{V}^{n}$ and $\mathcal{W}$ such that the following almost surely holds:
$$\left( \mathcal{V}^{n}, \hat{\mathcal{S}}_{\infty}^{n} \right) \overset{a.s.}{\to} \left( \mathcal{W}, \hat{\mathcal{S}}_{\infty} \right) \text{ in } \left( \mathcal{H} \times \Pi^{\mathbb{N}} , d_{\mathcal{H}} \times \hat{\mu}_{\infty} \right) \text{ as } n \to \infty$$
\begin{proof} Take the coupling of $\mathcal{V}^{n}, \mathcal{W}$ given in Theorem \ref{T5.3.2}. Fix any itinerary $\Xi_{k} = (x,s,\overline{\sigma}, \overline{\eta})$. For $0 \leq j \leq k$ let $\hat{\Xi}_{j}$ be a length $j$ restriction of $\Xi_{k}$. For convenience of notation, let $s=\sigma_{0}$ and $\infty=\sigma_{k+1}$. Note that by Theorem \ref{T4.0.2} we have:
$$\limsup_{n \to \infty} G_{\hat{\Xi}_{j}}^{n}(t) \leq \mathcal{G}_{\hat{\Xi}_{j}}(t) \text{ for } t > \sigma_{j}, \eta_{j} = 1 \text{ or } j=0$$
\begin{equation}\label{E7.3.2.1}
\limsup_{n \to \infty} -G_{\hat{\Xi}_{j}}^{n}(t) \leq -\mathcal{G}_{\hat{\Xi}_{j}}(t) \text{ for } t > \sigma_{j}, \eta_{j} = -1
\tag{7.5.1}
\end{equation}
Define paths $G_{\hat{\Xi}_{j}}^{\star, n}, \mathcal{G}_{\hat{\Xi}_{j}}^{\star}$ as paths from $[\sigma_{j}, \infty) \to \mathbb{R}$ in the following way:
$$G_{\hat{\Xi}_{j}}^{\star, n}(t) := \begin{cases} G_{\hat{\Xi}_{j}}^{n}(\sigma_{j}) & t \in [s, \sigma_{j}] \\ G_{\hat{\Xi}_{j}}^{n}(t) & t \in [\sigma_{j}, \sigma_{j+1}] \\  G_{\hat{\Xi}_{j}}^{n}(\sigma_{j+1}) & t > \sigma_{j+1} \end{cases}$$
$$\mathcal{G}_{\hat{\Xi}_{j}}^{\star}(t) := \begin{cases} \mathcal{G}_{\hat{\Xi}_{j}}(\sigma_{j}) & t \in [s, \sigma_{j}) \\ \mathcal{G}_{\hat{\Xi}_{j}}(t) & t \in [\sigma_{j}, \sigma_{j+1}) \\  \mathcal{G}_{\hat{\Xi}_{j}}(\sigma_{j+1}) & t > \sigma_{j+1} \end{cases}$$
By Theorem \ref{T7.3.1} we have:
$$G_{\Xi_{k}}^{n} \overset{d}{\to} \mathcal{G}_{\Xi_{k}} \text{ in } (\Pi, d) \implies$$
$$G_{\hat{\Xi}_{j}}^{\star, n} \overset{d}{\to} \mathcal{G}_{\hat{\Xi}_{j}}^{\star} \text{ in } (\Pi, d) \text{ for } 0 \leq j \leq k$$
Therefore by Lemma \ref{L9.0.7} and \eqref{E7.3.2.1} we have:
$$G_{\hat{\Xi}_{j}}^{\star, n} \overset{\mathbb{P}}{\to} \mathcal{G}_{\hat{\Xi}_{j}}^{\star} \text{ in } (\Pi, d) \text{ for } 0 \leq j \leq k$$
It follows that:
$$G_{\Xi_{k}}^{n} \overset{\mathbb{P}}{\to} \mathcal{G}_{\Xi_{k}} \text{ in } (\Pi, d)$$
Since the convergence is in probability, it can be simultaneously made for any finite number of itineraries $\Xi_{k}$. It follows by the definition of $\hat{\mu}_{\infty}$ that, under such a coupling:
$$\left( \mathcal{V}^{n}, \hat{\mathcal{S}}_{\infty}^{n} \right) \overset{\mathbb{P}}{\to} \left( \mathcal{W}, \hat{\mathcal{S}}_{\infty} \right) \text{ in } \left( \mathcal{H} \times  \Pi^{\mathbb{N}} , d_{\mathcal{H}} \times \hat{\mu}_{\infty} \right) \text{ as } n \to \infty$$
Applying the Skorokhod representation theorem then gives the desired result.
\end{proof}
\end{customthm}

\section{Appendix}
Here we include minor or tedious proofs whose arguments would disrupt the flow of the paper if included elsewhere.
\begin{customlemma}{8.1} \label{L9.0.1}  Let $\pi^{n}_{0}, ..., \pi^{n}_{k} \in \Pi$ and $\gamma_{0}, ..., \gamma_{k} \in \mathcal{W}$ and let $c > 0$. Suppose the following hold:
\begin{enumerate}
\item There exist jump times $w^{n}_{i}, t_{i}^{n}$ such that for all $n$ and for $1 \leq i \leq k$, the following hold:
$$\left| \pi_{i-1}^{n}(w_{i}^{n}) - \pi_{i}^{n}(t_{i}^{n}) \right| \leq c$$
$$w_{i}^{n}, t_{i}^{n} \in [-c, c]$$
\item There exist start times $\{ t^{n}_{0} \} \subset [-c, c]$ such that:
$$\lim_{n \to \infty} \pi^{n}_{0}(t^{n}_{0}) = x \text{ for some } x \in \mathbb{R}$$
\item $\pi_{i}^{n} \to \gamma_{i}$ in $(\Pi, d)$ as $n \to \infty$ for $0 \leq i \leq k$
\end{enumerate}
Then $\gamma_{i} \in \Pi^{\star}$ for all $i$.
\begin{proof} We proceed by induction on $i$.
\newline
\newline
\textbf{Base case:}
\newline
\newline
By the convergence of $\pi_{0}^{n} \to \gamma_{0}$ we have:
$$\inf\limits_{t \in [-c, c]} \gamma_{0}(t) \leq \limsup\limits_{n \to \infty} \inf\limits_{t \in [-c, c]} \pi^{n}_{0}(t) \leq x < \infty$$
Hence $\gamma_{0}(s) < \infty$ for at least some $s \in [-c,c]$. Similarly:
$$\sup\limits_{t \in [-c, c]} \gamma_{0}(t) \geq \liminf\limits_{n \to \infty} \sup\limits_{t \in [-c, c]} \pi^{n}_{0}(t) \geq x > -\infty$$
By the characterization given by \cite{FINR}, all paths in $\mathcal{W}$ are either identically $\infty$, identically $-\infty$, or elements of $\Pi^{\star}$. The above shows that $\gamma_{0}$ is neither identically $\infty$ nor identically $-\infty$. Hence $\gamma_{0} \in \Pi^{\star}$.
\newline
\newline
\textbf{Induction:} Suppose $\gamma_{i} \in \Pi^{\star}$ for some $i<k$. By way of contradiction, suppose further that $\gamma_{i+1} \notin \Pi^{\star}$. Then it is the case that either $\gamma_{i+1} \equiv \infty$ or $\gamma_{i+1} \equiv - \infty$. Without loss of generality, we address the case where $\gamma_{i+1} \equiv \infty$. Since $\pi_{i+1}^{n} \to \gamma_{i+1}$ in $(\Pi, d)$ we would then have:
$$\liminf_{n \to \infty} \inf\limits_{t \in [-c, c]} \pi_{i+1}^{n}(t) = \infty$$
But since $w_{i+1}^{n}, t_{i+1}^{n} \in [-c, c]$ are such that $\left| \pi_{i}^{n}(w_{i+1}^{n}) - \pi_{i+1}^{n}(t_{i+1}^{n}) \right| \leq c$, this would imply that:
$$\liminf_{n \to \infty} \sup\limits_{t \in [-c, c]} \pi^{n}_{i}(t) = \infty$$
Since $\pi_{i}^{n} \to \gamma_{i}$ in $(\Pi, d)$ it would follow that:
$$\sup\limits_{t \in [-c, c]} \gamma_{i}(t) = \infty$$
Contradicting the fact that $\gamma_{i} \in \Pi^{\star}$. Hence $\gamma_{i+1} \not\equiv \infty$. A symmetric argument shows $\gamma_{i+1} \not\equiv -\infty$. Hence it must be the case that $\gamma_{i+1} \in \Pi^{\star}$.
\end{proof}
\end{customlemma}
\text{ }
\newline

\begin{customlemma}{8.2} \label{L9.0.2}
Let $\mathcal{C} \subset \mathbb{R}$ be compact. Suppose $\pi^{n}$ converges to $\pi$ in $(\Pi, d)$ where $\pi \in \Pi^{\star}$. Then $\hat{\pi}^{n}$ converges uniformly to $\hat{\pi}$ along $\mathcal{C}$.
\begin{proof} Take $M>0$ such that $\mathcal{C} \subset [-M, M]$. Since $\pi \in \Pi^{\star}$, it is the case that $\hat{\pi}$ is bounded along $[-M, M]$. Thus there exists $N, \alpha > 0$ such that $1 > N+\alpha$ and $\sup\limits_{-M \leq t \leq M} |\tanh(\hat{\pi}(t))| \leq N$.
\newline
\newline
Since $[-(N+\alpha), N+\alpha]$ is compact, $\tanh^{-1}|_{[-(N+\alpha), N+\alpha]}$ is uniformly continuous. Fix any $\epsilon > 0$. Find $\delta_{1} > 0$ such that for any $w,z \in [-(N+\alpha), N+\alpha]$ the following holds:
$$|w-z| < \delta_{1} \implies |\tanh^{-1}(w)-\tanh^{-1}(z)| < \epsilon$$
Define:
$$\delta_{2} := \frac{\delta_{1}}{1+M} \wedge \frac{\alpha}{1+M}$$
Find $n_{0}$ such that for $n>n_{0}$ it is the case that $d(\pi^{n}, \pi) < \delta_{2}$. Then for such $n$ we have:
$$\sup\limits_{-M \leq t \leq M} |\tanh(\hat{\pi}^{n}(t))-\tanh(\hat{\pi}(t))| < \alpha \implies \tanh(\hat{\pi}^{n}(t)) \in [-(N+\alpha), N+\alpha] \text{ } \forall t \in [-M, M]$$
$$\sup\limits_{-M \leq t \leq M} |\tanh(\hat{\pi}^{n}(t))-\tanh(\hat{\pi}(t))| < \delta_{1} \implies \sup\limits_{-M \leq t \leq M} |\hat{\pi}^{n}(t) - \hat{\pi}(t)| < \epsilon$$
Hence for any $\epsilon > 0$, there exists $n_{0}$ such that for $n>n_{0}$ it is the case that $|\hat{\pi}^{n}(t) - \hat{\pi}(t)| < \epsilon$ for all $t \in \mathcal{C}$ as desired.
\end{proof}
\end{customlemma}
\text{ }
\newline

\begin{customlemma}{8.3}\label{L9.0.3} Equip $(\Pi^{\star})^{2}$ with the standard product metric $d^{2}$. Then $\Lambda_{R}$ and $\Lambda_{L}$ are continuous maps from $((\Pi^{\star})^{2}, d^{2}) \to (\Pi^{\star}, d)$.
\begin{proof}
We prove the theorem only for $\Lambda_{R}$, as the argument for $\Lambda_{L}$ is identical. Fix any $((f_{1}, t_{1}), (f_{2}, t_{2})) \in (\Pi^{\star})^{2}$ and any $\epsilon >0$. Let $\tau_{1} := t_{1} \vee t_{2}$. Find $\delta_{1} > 0$ such that:
$$|s-\tau_{1} | <\delta_{1} \implies \left| \hat{f}_{1}(s) - \hat{f}_{1}(\tau_{1}) \right| \vee \left| \hat{f}_{2}(s) - \hat{f}_{2}(\tau_{1}) \right| < \epsilon$$
Such a $\delta_{1}$ exists since $f_{1}, f_{2} \in \Pi^{\star}$. Find $M > 0$ such that $\frac{2}{1+M} < \epsilon$ and such that $[\tau_{1}-\epsilon, \tau_{1}+\epsilon] \subseteq [-M, M]$. Find some $\delta_{2} > 0$ such that:
$$d((f_{1}, t_{1}), (g_{1}, w_{1})) \vee d((f_{2}, t_{2}), (g_{2}, w_{2})) < \delta_{2} \implies$$
$$\sup\limits_{-M \leq t \leq M} \left| \hat{f}_{1}(t) - \hat{g}_{1}(t) \right| \vee\sup\limits_{-M \leq t \leq M} \left| \hat{f}_{2}(t) - \hat{g}_{2}(t) \right| \vee |t_{1}-w_{1}| \vee |t_{2} - w_{2}| \leq \delta_{1} \wedge \epsilon$$
Such a $\delta_{2}$ exists as by Lemma \ref{L9.0.2} convergence of $g^{n}$ to $f$ in $\Pi^{\star}$ implies uniform convergence of $\hat{g}^{n}$ to $\hat{f}$ along any compact set. For such a $\delta_{2}$, a simple computation shows that:
$$d((f_{1}, t_{1}), (g_{1}, w_{1})) \vee d((f_{2}, t_{2}), (g_{2}, w_{2})) < \delta_{2} \implies$$
$$d(\Lambda_{R}(f_{1}, f_{2}), \Lambda_{R}(g_{1}, g_{2})) < 5 \epsilon$$
\end{proof}
\end{customlemma}
\text{ }
\newline

\begin{customlemma}{8.4}\label{L9.0.4} Let $X^{n}, X$ be real random variables satisfying the following conditions:
\begin{enumerate}
\item $\limsup_{n \to \infty} X^{n} \leq X$ almost surely
\item $X^{n} \overset{d}{\to} X$
\end{enumerate}
Then $X^{n} \overset{\mathbb{P}}{\to} X$.
\begin{proof} 
By the assumed convergence in distribution, we have $\lim\limits_{n \to \infty} \mathbb{E}(\tanh(X^{n})-\tanh(X)) = 0$. As $\tanh(X^{n}), \tanh(X)$ are bounded, we apply Fatou's lemma giving:
\begin{equation*}
\begin{split}
0 &= \lim\limits_{n \to \infty} \mathbb{E}(\tanh(X^{n})-\tanh(X)) \\
&= \limsup_{n \to \infty} \mathbb{E}(\tanh(X^{n})-\tanh(X)) \\
&\leq \mathbb{E}(\limsup_{n \to \infty}\tanh(X^{n}) - \tanh(X)) \leq 0
\end{split}
\end{equation*}
Since $\limsup\limits_{n \to \infty} \tanh(X^{n}) - \tanh(X)$ is almost surely non-positive and has expectation zero, we have:
$$\limsup_{n \to \infty} \tanh(X^{n}) - \tanh(X) \overset{a.s.}{=} 0$$
Hence:
$$\limsup_{n \to \infty} X^{n} \overset{a.s.}{=} X$$
Define $Z^{n} := \sup_{m \geq n} X^{m}$. The above result shows $Z^{n} \overset{a.s.}{\to} X$. By definition, $Z^{n} \geq X^{n}$. Combining these observations we have:
\begin{equation*}
\begin{split}
\lim\limits_{n \to \infty}\mathbb{E}( |Z^{n}-X^{n}|) &= \lim\limits_{n \to \infty}\mathbb{E}(Z^{n}-X^{n}) \\
&= \lim\limits_{n \to \infty} \mathbb{E}(Z^{n}) - \mathbb{E}(X^{n}) \\
&= \mathbb{E}(X)-\mathbb{E}(X)=0
\end{split}
\end{equation*}
Thus $\lim\limits_{n \to \infty}\mathbb{E}(|Z^{n}-X^{n}|) =0$, which implies: 
$$\lim\limits_{n \to \infty} \mathbb{P}(|Z^{n}-X^{n}| > \epsilon) = 0$$ For any $\epsilon > 0$.
\newline
\newline
Since $Z^{n} \overset{a.s.}{\to} X$, it is the case that $Z^{n} \overset{\mathbb{P}}{\to} X$, and thus, for any $\epsilon >0$:
$$\lim\limits_{n \to \infty} \mathbb{P}(|Z^{n}-X| > \epsilon) = 0$$ 
\newline
\newline
Then for any $\epsilon >0$ we have:
$$\mathbb{P}( |X^{n} -X| > 2\epsilon) \leq \mathbb{P}(|X^{n}-Z^{n}| > \epsilon)+\mathbb{P}(|Z^{n}-X| > \epsilon) \overset{n \to \infty}{\to} 0$$
So that $X^{n} \overset{\mathbb{P}}{\to} X$ as desired.
\end{proof}
\end{customlemma}
\text{ }
\newline

\begin{customlemma}{8.5} \label{L9.0.5} Let $\delta, M > 0$. Let $f, g \in \Pi^{\star}$. Then the following hold:
\begin{itemize}
\item $\mathcal{M}_{\delta, M}(\Lambda_{R}(f,g)) \leq 2\mathcal{M}_{\delta, M}(f) + \mathcal{M}_{\delta, M}(g)$
\item $\mathcal{M}_{\delta, M}(\Lambda_{L}(f,g)) \leq 2\mathcal{M}_{\delta, M}(f) + \mathcal{M}_{\delta, M}(g)$
\end{itemize}
\begin{proof} We prove only the first statement as the argument for the second is identical. Suppose that $\mathcal{M}_{\delta, M}(\Lambda_{R}(f, g)) > \beta$. Then there exist $x,y$ satisfying:
\begin{enumerate}
\item $|x-y|<\delta$
\item $-M \leq x \leq y \leq M$
\item $\left| \Lambda_{R}(f, g)(y) - \Lambda_{R}(f,g)(x) \right| > \beta$
\end{enumerate}
First suppose that $\Lambda_{R}(f,g)(t) = g(t)$ for some $t \in [x, y]$. Define:
$$\tau_{1} := \inf \{ t \in [x,y] : \Lambda_{R}(f,g)(t) = g(t) \}$$
$$\tau_{2} := \sup \{ t \in [x,y] : \Lambda_{R}(f,g)(t) = g(t) \}$$
Then we have:
\begin{equation*}
\begin{split}
\beta &< \left| \Lambda_{R}(f, g)(y) - \Lambda_{R}(f, g)(x) \right| \\
& \leq \left| \Lambda_{R}(f, g)(y) - \Lambda_{R}(f,g)(\tau_{2}) \right| 
+ \left| \Lambda_{R}(f, g)(\tau_{2}) - \Lambda_{R}(f,g)(\tau_{1}) \right| \\
&+ \left| \Lambda_{R}(f, g)(\tau_{1}) - \Lambda_{R}(f, g)(x) \right|
\end{split}
\end{equation*}
Since $\Lambda_{R}(f,g)$ is a right Skorokhod reflection of $\hat{f}$ off of $\hat{g}$, by our selection of $\tau_{1}, \tau_{2}$, we have: 
$$\Lambda_{R}(f, g)(t) > \hat{g}(t) \text{ for all } t \in (x, \tau_{1}) \cup (\tau_{2}, y)$$
Further, by the properties of the Skorokhod reflection, we know that $\Lambda_{R}(f,g)$ evolves as $\hat{f}$ along time intervals for which $\Lambda_{R}(f, g)(t) \neq \hat{g}(t)$, so that:
$$\left| \Lambda_{R}(f,g)(\tau_{1})-\Lambda_{R}(f,g)(x) \right| = \left| \hat{f}(\tau_{1})-\hat{f}(x) \right| \leq \mathcal{M}_{\delta, M}(f)$$
Similarly,
$$\left| \Lambda_{R}(f,g)(y)-\Lambda_{R}(f,g)(\tau_{2}) \right| \leq \mathcal{M}_{\delta, M}(f)$$
Since elements of $\Pi^{\star}$ are continuous, we also have:
$$\left| \Lambda_{R}(f,g)(\tau_{2})-\Lambda_{R}(f,g)(\tau_{1}) \right| = \left| g(\tau_{2})-g(\tau_{1}) \right| \leq \mathcal{M}_{\delta, M}(g)$$
Hence in this case we see that:
$$\beta < \mathcal{M}_{\delta, M}(g)+ 2\mathcal{M}_{\delta, M}(f)$$
\newline
If instead, $\Lambda_{R}(f,g)(t) > g(t)$ for all $t \in [x,y]$, again by the properties of the Skorokhod reflection, we have:
\begin{equation*}
\begin{split}
\left| \Lambda_{R}(f,g)(y)-\Lambda_{R}(f,g)(x) \right| &= \left| \hat{f}(y)-\hat{f}(x) \right| \\
&\leq \mathcal{M}_{\delta, M}(f) \\
&\leq \mathcal{M}_{\delta, M}(g)+2\mathcal{M}_{\delta, M}(f)
\end{split}
\end{equation*}
In either case the result holds.
\end{proof}
\end{customlemma}
\text{ }
\newline

\begin{customlemma}{8.6} \label{L9.0.6}
Fix any $M, \delta >0$, then $\mathcal{M}_{\delta, M}$ is a continuous map from $(\Pi^{\star}, d) \to \mathbb{R}$.
\begin{proof}
Fix any $f \in \Pi^{\star}, \epsilon > 0$. It suffices to show that if $f^{n} \to f$ it is the case that $\mathcal{M}_{\delta, M}(f^{n}) \to \mathcal{M}_{\delta, M}(f)$.
\newline
\newline
Suppose $f^{n} \to f$ in $(\Pi^{\star}, d)$. By Lemma \ref{L9.0.2} $\hat{f}^{n} \to \hat{f}$ uniformly along $[-M, M]$. Thus, by taking $n$ large enough such that $\sup\limits_{t \in [-M, M]} |\hat{f}^{n}(t) - \hat{f}(t)| < \epsilon$ we have:
$$\left| \mathcal{M}_{\delta, M}(f) - \mathcal{M}_{\delta, M}(f^{n}) \right| < 2 \epsilon$$
\end{proof}
\end{customlemma}
\text{ }
\newline

\begin{customlemma}{8.7} \label{L9.0.7}
Let $(f^{n}, t^{n}), (f, t)$ be $\Pi$ valued random variables. Suppose that the following hold:
\begin{enumerate}
\item Almost surely for all $s \in \mathbb{R}$, it is the case that $\limsup_{n \to \infty} \hat{f}^{n}(s) \leq \hat{f}(s)$
\item $(f^{n}, t^{n}) \overset{d}{\to} (f, t)$ as elements of $(\Pi, d)$.
\item $t^{n} \overset{\mathbb{P}}{\to} t$
\end{enumerate}
Then $(f^{n}, t^{n}) \overset{\mathbb{P}}{\to} (f, t)$ in $(\Pi, d)$.
\begin{proof} Define $L_{0} : \Pi \to \mathbb{R}$ as:
$$L_{0}(f, t) := \tanh(t)$$
Let $\{ x_{j} \}_{j \in \mathbb{N}}$ be an enumeration of $\mathbb{Q}$. Define $L_{j} : \Pi \to \mathbb{R}$ as:
$$L_{j}(f, t) := \tanh(\hat{f}(x_{j}))$$
First we claim that $L_{j}$ is continuous for each $j$. To see this, fix any $\epsilon > 0$ and $j \in \mathbb{Z}_{\geq 0}$. Let $$\delta := \begin{cases} \frac{\epsilon}{1+|x_{j}|} & \text{ for }j \neq 0 \\ \epsilon & \text{ for } j=0 \end{cases}$$ Then if $d((f_{1}, t_{1}), (f_{2}, t_{2})) < \delta$, by definition $|L_{j}(f_{1}, t_{1}) - L_{j}(f_{2}, t_{2})| < \epsilon$. Hence $L_{j}$ is continuous as desired.
\newline
\newline
Next we claim that the collection $\{ L_{j} \}_{j \in \mathbb{Z}_{\geq 0}}$ is point separating. Suppose $L_{j}( f_{1}, t_{1}) = L_{j}(f_{2}, t_{2})$ for all $j \in \mathbb{Z}_{\geq 0}$. We will show if this holds, then $(f_{1}, t_{1})=(f_{2}, t_{2})$ as elements of $(\Pi, d)$.
\newline
\newline
First we have $L_{0}(f_{1}, t_{1}) = L_{0}(f_{2}, t_{2})$, hence $t_{1}=t_{2}$. 
\newline
\newline
Next take any $x \in \mathbb{R}$. Note that by definition of $\Pi$
, $\phi(\hat{f}(t), t)$ is continuous, hence:
\begin{equation*}
\begin{split}
\phi(\hat{f}_{1}(x), x) &= \lim\limits_{q \in \mathbb{Q}, q \to x} \phi(\hat{f}_{1}(q), q) \\
&= \lim\limits_{q \in \mathbb{Q}, q \to x} \frac{L_{j(q)}(f_{1}, t_{1})}{1+|q|} 
\\&= \lim\limits_{q \in \mathbb{Q}, q \to x}  \frac{L_{j(q)}(f_{2}, t_{2})}{1+|q|} \\
&= \phi(\hat{f}_{2}(x), x)
\end{split}
\end{equation*}
Hence $\phi(\hat{f}_{1}(x), x) = \phi(\hat{f}_{2}(x), x)$ for all $x \in \mathbb{R}$, thus $d((f_{1}, t_{1}), (f_{2}, t_{2})) = 0$ so that $(f_{1}, t_{1})=(f_{2}, t_{2})$ as desired.
\newline
\newline
Lastly we claim that $L_{j}(f^{n}, t^{n}) \overset{\mathbb{P}}{\to} L_{j}(f, t_{0})$ for all $j \in \mathbb{Z}_{\geq 0}$. The statement is true by assumption for $L_{0}$. For $j > 0$, by the continuous mapping theorem, we have:
$$(f^{n}, t^{n}) \overset{d}{\to} (f,t) \implies L_{j}(f^{n}, t^{n}) \overset{d}{\to} L_{j}(f, t)$$
By the monotonicity of $\tanh$ and our first assumption we have:
$$\limsup\limits_{n \to \infty} L_{j}(f^{n}, t^{n}) \overset{a.s.}{\leq} L_{j}(f, t)$$
Then by Lemma \ref{L9.0.4}, $L_{j}(f^{n}, t^{n}) \overset{\mathbb{P}}{\to} L_{j}(f, t)$.
\newline
\newline
Since $(f^{n}, t^{n}) \overset{d}{\to} (f, t)$, the sequence $\{ (f^{n}, t^{n}) \}_{n \in \mathbb{N}}$ is tight. By \cite{FINR}, $(\Pi, d)$ is a complete separable metric space. The result therefore follows from Proposition 3.12 in \cite{hair}.

\end{proof}
\end{customlemma}

\section*{Acknowledgements} I would like to acknowledge support from the Natural Sciences and Engineering Research Council of Canada (NSERC) through a postgraduate scholarship (PGS-D). Thanks to B\'alint Vir\'ag for proposing the question addressed in this paper, as well as for his extensive assistance in the problem-solving process and his recommendations on the presentation of the manuscript. Thanks to B\'alint Vet\H{o} for providing many thoughtful comments on content of the paper and and identifying points that needed clarification.

\bibliographystyle{dcu}
\bibliography{projectbib}

\end{document}